\documentclass[sn-vancuover,Numbered]{sn-jnl}

\usepackage{graphicx}%
\usepackage{multirow}%
\usepackage{amsmath,amssymb,amsfonts}%
\usepackage{amsthm}%
\usepackage{mathrsfs}%
\usepackage[title]{appendix}%
\usepackage{xcolor}%
\usepackage{textcomp}%
\usepackage{manyfoot}%
\usepackage{booktabs}%
\usepackage{algorithm}%
\usepackage{algorithmicx}%
\usepackage{algpseudocode}%
\usepackage{listings}%

\usepackage{breakurl}
\usepackage{underscore}
\usepackage[english]{babel}
\usepackage[style=numeric,minbibnames=3,natbib=false]{biblatex}
\usepackage{amssymb}
\usepackage{amsthm}
\usepackage{amsmath}
\usepackage{textcmds}
\usepackage{listings}
\usepackage{comment}
\usepackage{tikz-cd}
\usepackage{geometry}

\newcommand{\R}{\mathbb{R}}
\newcommand{\C}{\mathbb{C}}

\newcommand{\N}{\mathbb{N}}

\newcommand{\Set}{\mathbf{Set}}

\newcommand{\inv}{^{-1}}

\newcommand{\cata}{\textnormal{cata}}
\newcommand{\falg}{F_\Sigma\textnormal{-\textbf{Alg}}}

\newcommand{\ddd}{\cdot\cdot\cdot}

\theoremstyle{thmstyleone}%
\newtheorem{theorem}{Theorem}[section]
\newtheorem*{theorem*}{Theorem}

\newtheorem{lemma}[theorem]{Lemma}
\newtheorem{prop}[theorem]{Proposition}

\theoremstyle{thmstyletwo}%

\theoremstyle{thmstylethree}%
\newtheorem{definition}[theorem]{Definition}

\theoremstyle{definition}
\newtheorem{example}[theorem]{Example}

\begin{document}

\title[Algebraic Dynamical Systems in Machine Learning]{Algebraic Dynamical Systems in Machine Learning}

\author*[1,2]{\fnm{Iolo} \sur{Jones}}\email{iolo.j.jones@durham.ac.uk}
\author[2]{\fnm{Jerry} \sur{Swan}}\email{jerry@hylomorph-solutions.com}
\author[1]{\fnm{Jeffrey} \sur{Giansiracusa}}\email{jeffrey.giansiracusa@durham.ac.uk}

\affil[1]{\orgname{Durham University}, \orgaddress{\city{Durham}, \country{UK}}}
\affil[2]{\orgname{Hylomorph Solutions}, \orgaddress{\city{Glasgow}, \country{UK}}}

\abstract{We introduce an algebraic analogue of dynamical systems, based on term rewriting. We show that a recursive function applied to the output of an iterated rewriting system defines a formal class of models into which all the main architectures for dynamic machine learning models (including recurrent neural networks, graph neural networks, and diffusion models) can be embedded. Considered in category theory, we also show that these algebraic models are a natural language for describing the compositionality of dynamic models. Furthermore, we propose that these models provide a template for the generalisation of the above dynamic models to learning problems on structured or non-numerical data, including \q{hybrid symbolic-numeric} models.}




\maketitle

\section{Introduction}

The relationship between the structure of a model and its observable behaviour is central to many areas of applied mathematics. In linguistics and computer science, there are corresponding notions of the syntax and semantics of an expression or program. In machine learning, behaviour is determined via learned parameters while the structure of a model is typically considered to be prescribed by \emph{hyperparameters}. This is often categorified in the context of functional programming, where programs are viewed as maps in a category of data types. Here the syntax is specified by an \emph{algebraic data type}, on which, in a general setting, recursively-defined functions determine the semantics \cite{603bccba60774b4b92054ceb3651b481}. These perspectives can be combined in formal machine learning theory, where the categorical perspective forms a basis for describing \q{compositionality}: the properties of a model's components which are preserved under composition. Compositional modelling provides vital support for safe and causal inference, where unconstrained neural approaches are known to be lacking \cite{damour2020underspecification}.

In this paper, we develop this theory to include the increasingly popular class of models based on dynamical systems. We describe these models via universal algebra \cite{BurrisSankappanavar1981} and category theory and show that term rewriting systems \cite{baadernipkow1998} are the exact algebraic analogue of dynamical systems, but also explicitly encode the structure of the model in their expression. The rewrite rule corresponds to this syntax or structure, while the semantics are specified by a recursive function on that algebraic structure. The categorical setting also allows us to talk, in full generality, about which properties of models are preserved under composition. We use this to prove that rewriting models are naturally compositional, in the sense that the corresponding dynamical system will lift to any category with the appropriate structure.

\subsection{Structural constraints in machine learning.}
A proper appreciation of the role played by structural constraints requires a brief summary of the history of Artificial Intelligence (AI). Ever since its inception \cite{mccarthy-proposal-for-dartmouth-1955}, the field of artificial intelligence has been split between ostensibly competing \emph{symbolic} and \emph{connectionist} perspectives. The symbolic approach was initially favoured, exemplified by so-called `Good Old Fashioned AI' (GOFAI) \cite{russel2010}, that typically attempts to model the world in terms of rules which manipulate opaque  symbols via formalisms such as predicate calculus. Such approaches were ill-suited for modelling the noisy real world and suffered from the `knowledge elicitation bottleneck' in which domain experts were unable to provide adequate domain models at the desired high-level of representation. By the early 1990s, it became widely accepted that GOFAI had failed.

In contrast, the connectionist approach seeks to model the world numerically, notably via the universal function approximation properties of neural networks \cite{citeulike:3561150}. The approach has benefited from successive innovations, in particular the development in the late 1980s of the backpropagation algorithm \cite{rumelhart1986learning} for training multi-layer networks and the emergence of greater computing power in the 2000s. The currently dominant approach of supervised machine learning addresses the knowledge elicitation bottleneck by learning from a labelled training set. However, the underlying statistical mechanics of the learning mechanism mean that there is a tendency to model correlation rather than cause \cite{Scholkopf2022}. Historically the main challenge in machine learning has been to create models that are scalable, train well, and can approximate a suitable class of functions for the task at hand. There has been significant progress on these fronts in the past decades
but the methods used are generally \emph{black box} in that their behaviour cannot be explained and/or guaranteed. This has led to an increasing emphasis on developing models with structural constraints that guarantee particular properties, or make the model's behaviour easier to understand.

Structural constraints can bring several important benefits, including:
\begin{enumerate}
    \item Imposing domain-specific constraints on the model. This can include symmetries like translation equivariance, physical constraints when modelling a known physical system, or safety constraints that limit the behaviour of the model. 
    In contrast to \emph{black box} models, such constraints make the model's behaviour more interpretable/auditable.
    \item Making problems well posed and learnable. Fitting a model to data requires a constrained form of model to be well defined, and tighter constraints will require less data to train.
    \item Improving robustness by passing through dimensional bottlenecks. Data usually contain noise in the ambient space which can be reduced by a \emph{learned} compression into a smaller space where the geometry is more likely to be causal.
    \item Counteracting the \q{curse of dimensionality}. Sparse data in high dimensions will become more dense when reduced in dimension, allowing effective resampling from their distribution or interpolation between points.
\end{enumerate}
As we explore in the next section, there has been a general trend towards models that, for all the reasons above, are more tightly structurally constrained. Since these models are constructed via the composition of smaller parts, these structural properties can all be expressed as algebraic relations. These relations can encode things like dimension or parameter-sharing between functions and are also required to prove that the model has given properties, such as being equivariant under some group action. Group equivariance is the particular focus of \emph{geometric deep learning} \cite{bronstein2021geometric}, which aims to derive machine learning models that are equivariant under actions on the input data.

\subsection{Non-numerical data and structured models}

While most of the main machine learning models work on vector-valued data, many problems are best described with non-numerical or mixed data, such as graphs or discretely labelled data. In addition to the descriptive role that the algebraic perspective affords for machine learning theory, an algebraic representation can directly express \emph{hierarchically-structured} models, in which nodes can be labelled with arbitrary (i.e.\ symbolic and/or numeric) contextual data. 

It is widely acknowledged that deep learning has difficulty in generalizing beyond the training set \cite{damour2020underspecification} and the absence of first-class hierarchical structure in Deep Learning representations has been conjectured to be one of the main reasons behind this \cite{theroadtogeneralintelligence}, as well as leaving DL vulnerable to adversarial attack \cite{DBLP:journals/corr/GoodfellowSS14}. Of particular interest is the associated ability to represent hierarchical structure without the need for encoding and decoding steps, since these have an attendant prospect of ignoring important structural constraints. Cognitive linguists have long argued that a key cognitive capacity is the ability to represent and transform \emph{recursive} structure \cite{Marcus2001-MARTAM-10}, since such representations allow infinitely many propositions to be finitely described. Such induced structure can then be used as a base substrate for the computational implementation of cognitive mechanisms such as abstraction and analogy \cite{theroadtogeneralintelligence}.

\subsection{Universal algebra and compositionality}

We use universal algebra \cite{BurrisSankappanavar1981} as a language for describing the structure of machine learning models. We will specifically work with free term algebras to describe these models, which are algebraic constructions consisting of all the possible combinations of \emph{terms} generated by a predefined set of constants and functions. This set is called the signature of the algebra, and its elements are terms such as $x, y, f(\cdot), g(\cdot,\cdot,\cdot)$, ... which can be used to form combinations like $g(x,y,f(y))$. We can therefore separate this purely symbolic term from its evaluation, which we obtain by substituting particular values in for $x$ and $y$ in some set $X$ and particular functions $f: X \rightarrow X$ and $g: X^3 \rightarrow X$. The evaluation of terms is therefore done recursively, and so is naturally viewed in the context of functional programming, where the term algebra and set $X$ are both algebraic data types and the \emph{evaluation} is a catamorphism between them. We will treat this theory in the language of category theory, which formalises the idea that certain properties of the functions constructing the model can be preserved under composition.

In this framework, the internal structure of the model is made explicit in the term expression. This allows us to describe structural constraints in terms of the signature, and the form of the expressions we generate with it. On the other hand, by treating these algebras within the wider context of category theory, we can describe the compositional properties of models with the categories in which they lie.

\subsection{Dynamical systems and term rewriting}

While this algebraic approach seems to richly describe the structure and compositionality of static functions $f: X^n \rightarrow X^m$, an important class of machine learning models are the dynamical systems of the form
$$
G: X^n \rightarrow X^n \qquad f:X^n \rightarrow X^m
$$
leading to the dynamical process:
$$
\begin{tikzcd}
& x_0
& x_1
& x_2
&\\
& y_0 \arrow[u, "f"] \arrow[r, "G"]
& y_1 \arrow[u, "f"] \arrow[r, "G"]
& y_2 \arrow[u, "f"] \arrow[r, "G"]
& {}
\end{tikzcd}
$$
This class of models includes recurrence relations, recurrent neural networks, message-passing models on graphs, and diffusion models, but there is no clear notion of algebraic dynamics with which to fit them into this algebraic framework. This means there is no algebraic description with which to describe the structure and structural constraints on the dynamics of a model. It also means we cannot apply category theory to describe a model's compositionality, which is a particularly important consideration in the dynamic case as, if a property of the dynamics only holds approximately, the error will compound over time.

We will address this here, and show that in fact the notion of dynamical system is perfectly captured by the algebraic notion of \emph{term rewriting}. A rewrite rule is a relation that rearranges the subterms of a term in some specified manner, such as $f(x, f(y,z)) \mapsto f(f(x,y),z)$ which defines associativity.

We use this to construct a purely algebraic class of models called \emph{rewriting models} using the theory of term rewriting and functional programming. In our case, the rewrite rule will define the inductive composition of terms which, when evaluated recursively, will generate a dynamical system. Our main result is to prove that these rewriting models coincide with the usual definition of dynamical systems when the algebraic expressions are evaluated recursively: the rewrite rule thereby precisely captures the \textit{syntax} of the dynamics, onto which the learning process imposes a \textit{semantics} as a recursive function. In this sense, these models encode both syntax and semantics simultaneously, and so provide a unified language for talking about both. The informal version of the theorem is as follows.

\begin{theorem*}[Algebraic dynamical systems, Theorem \ref{main_thm}]
The class of dynamical systems embeds in the class of rewriting models, such that every rewriting model projects onto a dynamical system.
\end{theorem*}

In other words, we show that term rewriting is the correct notion of \textit{algebraic dynamical system}. These rewriting models are precisely the same as dynamical systems on the level of output, but they allow for a fully algebraic representation of the structural evolution over time. As we show in the next section, this is an incredibly important feature in practice, but currently lacks a formal representation. Furthermore, this description can be used prescriptively to impose structural and relational constraints in a compositional way. We show this in the following theorem, also given informally.

\begin{theorem*}[Compositionality, Theorem \ref{compositionality}]
If the sets and functions used in a rewriting model are all in a category $C$, where $C$ has the same coproduct as $\Set$, then the rewriting model is equivalent to one in $C$.
\end{theorem*}

This can be interpreted as a universality property, which states that the algebraic dynamics given by rewriting describe dynamical systems in every category $C$ with the appropriate coproduct. This includes all the main categories of interest in machine learning, and so proves that rewriting models are a fully compositional construction.

In the context of hybrid machine learning, the algebraic setting we use in this paper is type-agnostic, since the algebraic description of a model is treated independently of the data type of the output. By developing a general framework for algebraic dynamical systems, we propose that the rewriting models used in this paper provide a template for the future development of hybrid models. These two theorems provide a theoretical justification for this, because they show that rewriting models are indeed the correct symbolic analogue of dynamical systems, and that the important properties of the type (whatever it turns out to be) are preserved.
 
\section{Algebraic Constraints in Machine Learning}

In this section we will survey some different uses of algebraic constraints in popular machine learning models. By an \emph{algebraic} constraint we mean a rule that forces a given function in a model to take a specified form, such as passing through a bottleneck or sharing weights with other functions. We will emphasise that imposing algebraic constraints is the same as controlling the structure of the model, and that desirable structural properties of models can be expressed and specified algebraically.

\subsection{Static models}

We will first discuss \emph{static} machine learning (ML) models that do not evolve over time with the iterated application of some function.

\begin{example}[Neural networks]
Neural networks are a standard tool in machine learning tool for function approximation from data. They are composed of linear functions $f_i: R^{n_i} \rightarrow R^{n_{i+1}}$ called \emph{layers} with nonlinear \emph{activation} functions applied componentwise to all the layers
$$
\R^{n_{in}} \rightarrow \R^{n_{1}} \rightarrow \ddd \rightarrow \R^{n_{out}}
$$
They are ubiquitous because they train well and are universal approximators (in the sense that they are dense in $L^2$ \cite{SCARSELLI199815}), but are not \textit{compositional} in that they do not have many guaranteed properties besides continuity. In many image and temporal applications we would like the output to be space- or time-equivariant, which is achieved with a convolutional neural network in which the layer functions $f_i$ are the convolution of the input with a learned kernel. This also massively reduces the number of parameters needed per layer, something called \emph{parameter sharing} in which the small number of kernel parameters is unfolded to a whole layer.

So among the family of neural networks, the ones with important compositional properties can be specified by this particular formula for convolution. The structure we want is therefore induced algebraically, and the algebraic form also controls things like the size and shape of the kernel.
\end{example}

\begin{example}[Low-rank matrix completion]
The matrix completion problem is fundamental in many applied ML domains \cite{7426724}. Given an incomplete matrix $M$, the problem is to find a complete matrix $M'$ that approximately agrees with $M$ where $M$ is sampled. This is a common problem in recommendation systems where each column is a product and each row is a user, and the matrix contains that user's rating or opinion of that product. Most users will not have used most products, so $M$ is mostly empty, but if it can be completed then the entries of $M'$ can be used to suggest products to a user based on their previous preferences.

The most common approach is to approximate $M$ by a low-rank matrix $M'$ so that, if $M$ is $n \times m$, $M' = AB$ where $A$ is $n \times k$ and $B$ is $k \times m$. If $k$ is taken to be small (less than the number of completed entries of $M$), then finding $A$ and $B$ becomes a well-posed optimisation problem. Furthermore, by tuning $k$ we can control the properties of the output: smaller $k$ will make $M'$ smoother and more robust while larger $k$ will make it more sensitive and precise. In this example the form of the solution is straightforward, but these desirable qualities are entirely structural, and are controlled by the algebraic properties of $M'$ (namely its rank).
\end{example}

\begin{example}[Autoencoders]
There are many \emph{generative} problems in ML where we would like a model to sample new data from a given probability distribution. These generally take the form of an autoencoder \cite{math11081777}, in which data in a large dimensional space are \emph{encoded} in a smaller dimensional space by a function $e: \R^n \rightarrow \R^k$, and then \emph{decoded} by a function $d: \R^k \rightarrow \R^n$. The model is then trained to minimise $||d \circ e - id||$, so that the encoding has minimal loss but is constrained by the bottleneck dimension $k$. Unlike the matrix completion example, there is not a natural choice for the forms of $e$ and $d$, although the bottleneck $k$ plays the same role of controlling precision versus robustness. Again, this $k$ is a structural constraint that is implemented through the algebraic form of the model.

The bottleneck is also essential for using (variational) autoencoders as generative models, where we want to sample data from a distribution by sampling in the compressed latent space. The data in the input space $\R^n$ often suffer from the \emph{curse of dimensionality}: that the number of data needed to densely sample a region of space increases exponentially with the dimension. Statistical modelling of data in this space is then usually impossible, but the sparse data in $\R^n$ can be made sufficiently dense in $\R^k$ for small enough $k$.
\end{example}

\begin{example}[Graph models]
Graph neural networks (GNNs) \cite{zhou2021graph} are a model for ML on data in the form of a labelled graph, where each vertex has a vector label $m_v$ with information about that vertex, and edges are sometimes labelled as well. This can be used for vertex segmentation or regression, as well as whole-graph classification. GNNs are static models in the above sense, where the data on each vertex is transformed by a function which depends both on the data at that vertex and the data of its neighbours (and the edges connecting the neighbours, if they are labelled). This constitutes one \emph{layer} of the network, and a GNN usually comprises a few layers.

In this sense, GNNs make use of the graph structure via the neighbourhood aggregation, and so constrain the algebraic form of the model to one that encodes locality. By placing a topological requirement on the algebraic expressions in the model, the GNN becomes sensitive to that extra structure and does not simply view the vertices as a set.

\end{example}

\subsection{Dynamic models}

We now consider dynamic ML models which are the main focus of this paper. By \emph{dynamic} we mean a model of which a large component is the iterated application of a learned function, such as a recurrence relation, recurrent neural network (RNN), or attention transformer. As discussed above, the algebraic properties of a static model specify its structure, and that structure is often essential for a variety of objectives. Historically the dominant trend in dynamical models was to use fairly unstructured dynamics, with minimal algebraic constraints, like general message-passing neural networks for graphs and RNNs for time series. More recently huge progress has been made by replacing these models with simpler but more structured alternatives with tighter algebraic constraints, such as convolutional networks \cite{8286426}, attention transformers \cite{NEURIPS20199d63484a} and diffusion models for graphs \cite{chamberlain2021grand}, attention head networks for time series \cite{NIPS20173f5ee243}, and the emergence of stable diffusion for generative modelling \cite{cao2022survey}. In all these cases the power of the model is derived from structural constraints on the dynamical interaction taking place, suggesting that an algebraic approach will effectively describe the properties of these models as well.

\begin{example}[Time series models]

The canonical example of a dynamic ML model is an RNN \cite{lipton2015critical}, or variations thereof such as \emph{long short-term memory} networks \cite{8737887}. In these models, there is a \emph{hidden} dynamical system $\R^k \rightarrow \R^k$ which is updated at each time step by a function that takes in the state vector from the previous time step and also the network's input, if it is made available. The \emph{update} function is fixed, so the algebraic form of the model reuses this base function at each timestep, leading to a family of models parametrised by time. This algebraic property guarantees that the model defines a dynamical system, and is effectively translation equivariant.
\end{example}

\begin{example}[Diffusion models]
Many popular ML models are trained by continuously optimising some objective function, and so must be differentiable in their parameters. They are typically also differentiable in their input (i.e.\ the models are differentiable functions), and so, once trained, can be used to optimise some \textit{input} to have a particular output value. A recent and very prominent extension of this idea is a \emph{diffusion model}, which takes random input and optimises the likelihood of that input given by some target probability distribution. This allows the generation of sample data from that distribution by starting with a random input point and then \emph{diffusing} it in this way towards the distribution. A more advanced variant uses Bayesian optimisation with a pre-trained classification model to generate prompt-specific diffusions to the distribution \cite{cao2022survey}.

Viewed as a discrete diffusion process, diffusion models are also dynamic models in the sense that the diffusion operator is iteratively applied to the input. The iteration of this operator thereby induces a family of models indexed by diffusion time, where the output of the model is taken to be the limit of this process in some sense. The theoretical guarantees for the diffusion follow from Langevin dynamics, which are encoded in the algebraic form of the diffusion operator. When Bayes' theorem is used for prompt-based diffusion, this is also specified algebraically. As with RNNs, this allows an algebraic description of the model's individual static functions, but the model itself, viewed as a family of $n$-step iterations, does not admit a unified algebraic description.
\end{example}

While we can use the algebraic theory for static models to describe the functions involved in these dynamic examples, there is no unified algebraic framework for expressing the model as a whole. In other words, there is no corresponding \emph{algebraic dynamics} that describes the iterated application of the component functions that comprise all these examples. We will address that question in this paper.

This is also an important issue in the development of hybrid dynamic models, which need to process timeseries of symbolic or mixed data types. For example, one particular open problem in contemporary robotics is how to extend 
`simultaneous localization and mapping' (SLAM \cite{doi:10.1177/027836498600500404}) algorithms (which incrementally construct scene representations in real-time) to incorporate semantic constraints \cite{DBLP:journals/corr/abs-1803-11288,hughes2023foundations}.
In the proposed model, this is possible by propagating successive parse-trees as first-class objects through the learning process for the timeseries.
The output of a parse is a structured representation of a scene and (in all but pathological situations) there will likely be significant structural continuity between frames. With a hybrid representation, this commonality can be exploited by the learning process to ensure spatio-temporal coherence and more efficient inference. 

\section{Related work}
As first introduced in Risi Kondor's thesis \cite{kondor2008group}, the term `algebraic machine learning' refers to harmonic approaches to representation-theoretic constraints; Swan \cite{7969355} relatedly used the group-theoretic Fourier transform to perform continuous (i.e.\ real-valued) learning of heuristics for a discrete (i.e.\ permutation) problem. A wide perspective on constraint representation for ML is proposed by Bronstein et al's 
`Erlangen programme' for Geometric Deep Learning via the unifying perspective of group theory \cite{bronstein2021geometric}. 

Following an initial miscellany of influential publications (e.g.\ \cite{10.1145/3110271,Elliott-2018-ad-icfp,fong2019backprop}), there is nascent but increasingly-convergent interest in the application of category theory to machine learning, with Shiebler et al providing a recent survey \cite{DBLP:journals/corr/abs-2106-07032} of this rapidly-growing area. An area of considerable activity is a categorical treatment of inference (e.g.\ \cite{smithe2020bayesian,cruttwell2021categorical,DBLP:journals/corr/abs-1911-12904}), in particular affording a unified perspective via (dependent) optics/polynomial functors (e.g.\ \cite{fong2019lenses,theroadtogeneralintelligence,2005.01894,Spivak2022,MyersCategoricalSystemsTheory}).

Regarding the desire to maintain a compositional mapping between syntax and semantics, previous work by Bloom et al 
\cite{DBLP:journals/acs/BloomSW96} (motivated predominantly by the ability to reason about concurrent systems) provides a functorial mapping from datatypes (the syntax) to behaviour (as represented by e.g.\ finite state automata). A compositional theory for the operational semantics of generalised Petri nets is proposed by Master \cite{master2021composing}. With particular regard to the role of algebraic data types as a syntactic form for structural learning, Swan \cite{10.1007/s10710-019-09347-3} uses symbolic regression to learn recursive functions of algebraic data types (ADTs) and subsequently proposes  \cite{theroadtogeneralintelligence} that the ADT structure itself should also be learned, as a basis for `necessary and sufficient' causal structure.


\section{Discrete dynamical systems}

The class of models we will seek to describe algebraically are known as discrete dynamical systems. We will regard these as models for time series, although they can have different interpretations, such as diffusion processes as described above. The term \emph{discrete dynamical system} is often also used more specifically for the sub-class which, in this paper, we will term \textit{recurrence relations} \cite{galor2005discrete}. These are the sequences $(s_n)$ which satisfy
$$
s_n = f(s_{n-1},...,s_{n-d})
$$
for some $d$ (which we call the \textit{depth} of the relation) and function $f: X^d \rightarrow X$, for all $n \geq d$. The sequence is then specified by $f$ and the $d$ initial conditions $s_0,...,s_{d-1} \in X$.

We will briefly discuss some examples and properties of recurrence relations, and how they fit within the broader class of recurrent systems.

\subsection{Recurrence relations}

A simple example of a recurrence relation is the Fibonacci sequence given by $s_n = s_{n-1} + s_{n-2}$, and linear recurrence relations in general are a surprisingly expressive class of model (i.e. where the function $f$ is linear). They include, for example, the popular ARIMA model for time series (albeit without the error quantification for which ARIMA is also used \cite{arima}).

\begin{example}
A sequence of the form
$$
s_n = \sum_{i=1}^m \big(a_i \sin(c_i n + s_i) + b_i \cos(c_i n + s_i)\big)
$$
can be expressed as a depth-$2m$ linear recurrence relation. Notice that
\begin{equation*}
    \begin{split}
        s_{n-k} &= \sum_{i=1}^m a_i \big(\sin(c_i n + s_i)\cos(c_i k) - \cos(c_i n + s_i)\sin(c_i k) \big)\\
        &\qquad + \sum_{i=1}^m b_i \big(\cos(c_i n + s_i)\cos(c_i k) - \sin(c_i n + s_i)\sin(c_i k) \big) \\
        &= \sum_{i=1}^m \big(a_i\cos(c_i k) - b_i\sin(c_i k)\big) \sin(c_i n + s_i)\\
        &\qquad + \sum_{i=1}^m \big(b_i\cos(c_i k) - a_i\sin(c_i k) \big)\cos(c_i n + s_i) \\
    \end{split}
\end{equation*}
for all $k \in\N$. Hence the set $\{s_{n-1},...,s_{n-2m}\}$ satisfies a system of linear equations in $2m$ variables, which can be solved to give each $\sin(c_i n + s_i)$ and $\cos(c_i n + s_i)$ in terms of the $s_i$. These can be substituted into the expression for $s_n$ to give a linear recurrence relation.
\end{example}

A broader class of sequences can be modelled by non-linear recurrence relations.

\begin{example}
Let $p$ be a degree-$k$ polynomial, and $s_n = p(n)$. Notice that we can find constants $b_0,...,b_k$ such that 
$$
b_0 + b_1p(n-1) + ... + b_kp(n-k) = n
$$
for all $n$, so define $f:\R^k\rightarrow\R$ by
$$
f(x_1,...,x_k) = b_0 + b_1x_1 + ... + b_kx_k.
$$
Then $s_n = p\circ f(s_{n-1},...,s_{n-k})$, with initial conditions $p(0),...,p(k)$.
\end{example}

\begin{example}
If $s_n$ is a strictly monotonic sequence in $\R$ then we can write $s_n = f(n)$ for some (non-unique) invertible function $f:\R \rightarrow \R$. So $s_n = f(f \inv (s_{n-1}) + 1)$ is a first order recurrence relation.
\end{example}

If the dynamical system is a recurrence relation, then the learning problem is greatly simplified. The structure of the model as $s_n = f(s_{n-1},...,s_{n-d})$ means that, for a given depth $d$, we can compute the time-delay embedding
$$
\{(s_n, s_{n-1},...,s_{n-d}) : n \in \N \}
$$
and \emph{learning} the recurrence relation $f$ becomes a regression problem on this point cloud.

\subsection{Dynamical Systems}

A common method in the analysis of dynamical systems is to exchange \emph{1-dimensional/ $n^{th}$ order} systems for \emph{$n$-dimensional/ first order} ones. For example, with a linear recurrence of order 2, we would write
$$
\begin{pmatrix}
s_n \\
s_{n-1}
\end{pmatrix}
=
\begin{pmatrix}
a_1 & a_2 \\
1 & 0
\end{pmatrix}
\begin{pmatrix}
s_{n-1} \\
s_{n-2}
\end{pmatrix}.
$$
Now we have a hidden \textit{first order} dynamical system which is \emph{interpreted} via the projection map $(x,y) \mapsto x$. This more general form of model is what we will call \emph{dynamical systems} in this paper. Rather than modelling a time series itself recurrently as above, there is a \emph{hidden} dynamical process in a latent space, along with a map from the system to the output space.

\begin{definition}
A \textbf{dynamical system} is a pair of sets $(X, Y)$ and maps $G: Y \rightarrow Y$ and $f: Y \rightarrow X$. We call $Y$ the \textbf{latent space (or set)} and $X$ the \textbf{output space (set)}. We call a dynamical system \textbf{cartesian} if $Y = X^n$ for some $n$.  
\end{definition}

A dynamical system produces a sequence in $X$ for each initial \emph{state} $y_0 \in Y$, by the iterated application of $G$ to $y_0$ followed by $f$.
$$
\begin{tikzcd}
& x_0
& x_1
& x_2
&\\
& y_0 \arrow[u, "f"] \arrow[r, "G"]
& y_1 \arrow[u, "f"] \arrow[r, "G"]
& y_2 \arrow[u, "f"] \arrow[r, "G"]
& {}
\end{tikzcd}
$$

\begin{example}
A recurrent neural network is a cartesian dynamical system in which $G$ and $f$ are neural networks, and usually $X = \R$.
\end{example}
\begin{example}
A Kalman filter's estimation process is a cartesian dynamical system in which $G$ and $f$ are linear, and $X = \R$ \cite{kalman1960}.
\end{example}
\begin{example}
A depth-$d$ recurrence relation is a cartesian dynamical system where $f$ is the projection $\pi_1$, and $G$ is given by
$$
(x_1,...,x_d) \mapsto (g(x_1,...,x_d), x_1,...,x_{d-1})
$$
\end{example}

\begin{example}
A (time-independent) message passing neural network is a dynamical system on the set of edges and vertices, but it can also be viewed as a cartesian dynamical system. Let each vertex have hidden state $h_v^t$ at time $t$, the edges have fixed hidden states $e_{vw}$, the messages be passed by
$$
m_v^{t+1} = \sum_{w \in N(v)} M(h^t_v, h^t_w, e_{vw})
$$
and states updated by
$$
h_v^{t+1} = U(h^t_v, m^{t+1}_v),
$$
with \emph{readout} function
$$
R(\{h_v^t : v \in G\}).
$$
If there are $n$ vertices with $d$-dimensional hidden states, then the recurrent system has hidden space $Y = \R^{nd}$, where each hidden state is of the form
$$
x_t = (h_{v_1}^t, ..., h_{v_n}^t).
$$
Define $m: Y \rightarrow \R^n$ by
$$
(h_1,...,h_n) \mapsto \big( \sum_{v_k \in N({v_1})} M(h_1, h_k, e_{v_1 v_k}), ..., \sum_{v_k \in N({v_n})} M(h_n, h_k, e_{v_n v_k}) \big),
$$
and $G: Y \rightarrow Y$ by 
$$
x = (h_1,...,h_n) \mapsto \big( U(h_1, m(x)),\ ...\ ,\ U(h_n, m(x)) \big).
$$
Then the dynamical system uses $G$ for state transitions and the output function $f:Y \rightarrow \R$ given by
$$
(h_1,...,h_n) \mapsto R(\{h_k : k=1,...,n\}).
$$
\end{example}

There are many other such examples of dynamical systems in machine learning, and the class is strictly greater than the recurrence relations. However, we can ask exactly how much bigger, or equivalently under what conditions a dynamical system can be expressed as a recurrence relation.

\begin{lemma}\label{rr_reduction}
Let $D$ be a cartesian dynamical system given by $G$ and $b$ and define $\phi: \R^d \rightarrow \R^d$ by
$$
x \mapsto (b(x), b(G(x)), ... , b(G^{d-1}(x))).
$$
Then $D$ is a depth-$d$ recurrence relation if $\phi$ is invertible.
\end{lemma}
\begin{proof}
If we define $\phi$ as above then
\begin{equation*}
    \begin{split}
        \phi(x_{k-d-1})
        &= (b(x_{k-d-1}),...,b(x_{k-2})) \\
        &= (s_{k-d},...,s_{k-1}),
    \end{split}
\end{equation*}
Where $s_n$ is the sequence induced by $D$, and $x_n$ is the sequence of hidden states. If $\phi$ is invertible we can recover $x_{k-d-1}$, so then
$$
s_k = b(G^{d}(\phi\inv(s_{k-d},...,s _{k-1}))),
$$
so we have a recurrence relation $b\circ G^{d}\circ\phi\inv$.
\end{proof}

In the case that the dynamical system is linear, we can directly apply this to \textit{generically} reduce the model to a recurrence relation.

\begin{prop}
If $K$ is a linear cartesian dynamical system where $G$ is linear with distinct eigenvalues, and $b$ acts non-trivially on the eigenvectors of $G$, then $K$ is a depth-$d$ linear recurrence relation, where $d$ is the dimension of the hidden space.
\end{prop}
\begin{proof}
Let $G : \R^d \rightarrow \R^d$ where $\R^d$ is the hidden state, and define $\phi:\R^d \rightarrow \R^d$ as the linear map
$$
x \mapsto (b(x), b(G(x)), ... , b(G^{d-1}(x))).
$$
Lemma \ref{rr_reduction} says that if $\phi$ is invertible then $K$ is a recurrence relation. Since $G$ has distinct eigenvalues $\lambda_j$, we can change the basis of $\C^d$ to an eigenbasis for $G$, where now $G$ is diagonal with distinct entries. So, writing $b = (b_1,...,b_d)$ in this new basis, where we assume all $b_i \neq 0$, we have
\begin{equation*}
    \begin{split}
        \det(\phi)
        &= \begin{vmatrix}
        b_1 & b_2 & \ddd & b_d \\
        b_1\lambda_1 & b_2\lambda_2 & \ddd & b_d\lambda_d \\
        \ddd & \ddd & \ddd & \ddd \\
        b_1\lambda_1^{d-1} & b_2\lambda_2^{d-1} & \ddd & b_d\lambda_d^{d-1} \\
        \end{vmatrix} \\
        &= \Big( \prod_{j=1}^d b_j \Big) \begin{vmatrix}
        1 & 1 & \ddd & 1 \\
        \lambda_1 & \lambda_2 & \ddd & \lambda_d \\
        \ddd & \ddd & \ddd & \ddd \\
        \lambda_1^{d-1} & \lambda_2^{d-1} & \ddd & \lambda_d^{d-1} \\
        \end{vmatrix} \\
        &= \Big( \prod_{j=1}^d b_j \Big) \Big( \prod_{1\leq i < j \leq d} (\lambda_j - \lambda_i) \Big),
    \end{split}
\end{equation*}
using the formula for the determinant of a Vandermonde matrix. This is non-zero under precisely the conditions we have assumed. So $\phi$ is invertible, and
$$
s_k = b \circ G^{d} \circ \phi\inv(s_{k-d},...,s_{k-1}),
$$
making $K$ a linear recurrence relation.
\end{proof}

This result can be seen as a special case of Takens' theorem \cite{takens2006detecting}, which says that the map $\phi$ is invertible, for some depth $d$, for a generic smooth dynamical system on $\R$. As such, we can view smooth dynamical systems as recurrence relations based on the principle in lemma \ref{rr_reduction}.

However, to say that a dynamical system \emph{can} be modelled as a recurrence relation does not mean that this is an efficient model for learning. Based on the principle of compositionality, in practice such models are often composed of more complex or robust parts, as described in the examples above. In the following section, we lay out our algebraic model for dynamical systems that offers a natural language for describing and reasoning about compositionality in time-series modelling.

\section{Universal Algebra and Rewriting}

We now introduce the proposed \emph{rewriting model} using ideas from functional programming. The model will comprise an iterated term rewriting process, in which a single term is iteratively rewritten to generate a sequence. This provides a purely algebraic representation of the structure of the model, which can then be interpreted by a recursive function from the term algebra into the output type. We will conclude that there is a constructive correspondence between the two families of model, so that, given a model in one family, we can construct an equivalent model in the other.

In this way we can identify the model architecture with two algebraic objects: the initial term at the start of the process, and the rewrite rule that tells us how to build up the dynamics at each time step. The recursively defined output function then \emph{interprets} the term in the output type. In other words, the rewrite rule describes the hyper-parameters of the model, including any algebraic constraints, and the catamorphism describes the parameters.

This will allow a unified algebraic description of the time-parametrised family of models and their constraints, and permits the same algebraic analysis and specification that exists for static models.

\subsection{Terms and Algebras}

The rewriting process is purely formal and will be defined on a term algebra. We now introduce some standard notions from universal algebra and term rewriting. More details can be found in \cite{baadernipkow1998}.

\begin{definition}[Signature]
A \textbf{signature} is a function
$\Sigma : \N \rightarrow \Set$. The set $\Sigma_n := \Sigma(n)$ is called the set of n-ary operators of $\Sigma$.
\end{definition}

\begin{definition}[$\Sigma$-terms]
Let $\Sigma$ be a signature and $V$ be a set of variables such that $\Sigma \cap V = \varnothing$. The set $T(\Sigma,V)$ of all $\Sigma$\textbf{-terms} over $V$ is defined inductively by
\begin{itemize}
    \item $V \subseteq T(\Sigma,V)$ (every variable is a term),
    \item $f(t_1, ..., t_n) \in T(\Sigma,V)$ for all $f \in \Sigma_n$ and $t_i \in T(\Sigma,V)$ for $i=1,...,n$ (closure under the application of function symbols).
\end{itemize}
\end{definition}

We can identify the terms with the set of planar trees whose vertices are labelled by elements of $\Sigma_n$ and leaves by elements of $V$. We can alternatively interpret this term set as the elements of a type, where the signature induces the constants and operators associated with the type. The variables will let us define equations and relations in the type with which we can define rewriting. For example, $\Sigma_0 = \{a,b\}$, $\Sigma_1 = \{f\}$, and $\Sigma_2 = \{g\}$ means that there are two constants, one unary function, and one binary function.

We will use the following standard definitions, although intuitively these correspond to operations on the tree representation of terms. Positions in the term index the vertices of the tree, a subterm is obtained by pruning the tree at a vertex, and a tree can be inserted into another at a specified vertex.

\begin{definition}
Let $\Sigma$ be a signature, $V$ be a set of variables disjoint from $\Sigma$, and $s,t \in T(\Sigma, V)$.
\begin{enumerate}
    \item The set of \textbf{positions} of the term $s$ is a set $Pos(s)$ of strings over the alphabet of positive integers, which is inductively defined as follows:
    \begin{itemize}
        \item If $s = v \in V$, then $Pos(s) := \{\epsilon\}$, where $\epsilon$ denotes the empty string.
        \item If $s = f(s_1,... ,s_n)$, then
        $$Pos(s) := \{\epsilon\} \cup \bigcup \{ip : p \in Pos(s_i)\}.$$
    \end{itemize}
    The position $\epsilon$ is called the \textbf{root position} of the term $s$, and the function or variable symbol at this position is called the \textbf{root symbol} of $s$. We denote by $l(p)$ the length of the string $p$.
    
    
    \item For $p \in Pos(s)$, the \textbf{subterm} of $s$ at position $p$, denoted by $s|_p$, is defined by induction on the length of $p$:
    \begin{itemize}
        \item $s|_\epsilon = s$,
        \item $f(s_1,...,s_n)|_{iq} = s_i|_q$.
    \end{itemize}

    \item For $p \in Pos(s)$, we denote by $s[t]_p$ the term that is obtained from $s$ by replacing the subterm at position $p$ by $t$, i.e.
    \begin{itemize}
        \item $s[t]_\epsilon = t$,
        \item $f(s_1,...,s_n)[t]_{iq} = f(s_1,...,s_i[t]_q,...,s_n)$.
    \end{itemize}
    
    \item We denote the set of variables occurring in $s$ by $Var(s)$, so
    $$
    Var(s) := \{v \in V : s|_p = v \textnormal{ for some } p \in Pos(s)\}.
    $$
    A term $t \in T(\Sigma, V)$ is called ground if $Var(t) = \varnothing$.
\end{enumerate}
\end{definition}

We can induce maps on the term algebra by defining them on the variables, and then insisting that these \emph{substitutions} commute with the operators.

\begin{definition}[Substitution]
Let $\Sigma$ be a signature and $V$ be a countably infinite set of variables. A \textbf{$T(\Sigma, V)$-substitution} (or just substitution), is a function $V \rightarrow T(\Sigma, V)$. The set of all $T(\Sigma, V)$-substitutions will be denoted by $Sub(T(\Sigma, V))$ or just $Sub$. Any $T(\Sigma, V)$-substitution $\sigma$ can be extended to a map $\hat{\sigma} : T(\Sigma, V) \rightarrow T(\Sigma, V)$ by setting $\hat{\sigma} = \sigma$ on $V$ and then inductively defining
$$\hat{\sigma}(f(s_1,... , s_n)) = f(\hat{\sigma}(s_1),... , \hat{\sigma}(s_n))$$
for all $s_i \in T(\Sigma, V)$. We say a term $t$ is an \textbf{instance} of a term $s$ if there exists a substitution $\sigma$ such that $\hat{\sigma}(s) = t$. We will usually suppress the hat in notation where there is no ambiguity.
\end{definition}

In functional programming, we can interpret the signature as generating an \emph{algebraic data type (ADT)} where the arity of the operators specify the \emph{shape} of the type in some sense. We now make this notion precise.

\begin{definition}[$\Sigma$-algebra]
If $\Sigma$ is a signature, a $\Sigma$\textbf{-algebra} is a set $X$ along with functions $x_\sigma: X^n \rightarrow X$ for all $\sigma \in \Sigma_n$ and all $n$. In particular, there are constants $x_\sigma \in X$ for each $\sigma \in \Sigma_0$.
\end{definition}

\begin{definition}[$\Sigma$-algebra Homomorphism]
If $\Sigma$ is a signature and $X$ and $Y$ are $\Sigma$-algebras, a $\Sigma$\textbf{-algebra homomorphism} is a function $h:X \rightarrow Y$ such that $h(x_\sigma(x_1,...,x_n)) = y_\sigma(h(x_1),...,h(x_n))$ for all $\sigma \in \Sigma_n$. In particular $h(x_\sigma) = y_\sigma$ for all $\sigma \in \Sigma_0$.
\end{definition}

So the unique extension of substitution (explained in their definition) is really saying that a function $f: V \rightarrow T(\Sigma, V)$ extends uniquely to a homomorphism $T(\Sigma, V) \rightarrow T(\Sigma, V)$. The $\Sigma$-algebras also form a category. The sets of $\Sigma$-terms $T(\Sigma, V)$ are $\Sigma$-algebras, where for each $f \in \Sigma_n$ we define $x_f: X^n \rightarrow X$ by
$$
(t_1,...,t_n) \mapsto f(t_1,...,t_n).
$$
It can be shown that the set of \textbf{ground} terms $T(\Sigma, \varnothing)$ is the initial object in this category, called the initial algebra. It is a subalgebra of every set of terms $T(\Sigma, V)$. From now on we will suppress the $\Sigma$ in the notation and denote, for fixed $\Sigma$, $T:= T(\Sigma, \varnothing)$ and $T(V):= T(\Sigma, V)$. We can equivalently describe the $\Sigma$-algebras as algebras of an endofunctor $F$ of $\Set$.

\begin{definition}[$F$-algebra]
If $C$ is a category, and $F:C\rightarrow C$ is an endofunctor of $C$, then an $F$\textbf{-algebra} is a tuple $ (A,\alpha)$, where $A$ is an object of $C$ and $\alpha $  is a $C$-morphism $F(A)\rightarrow A$. The $F$-algebras form a category where a map from $(A,\alpha)$ to $(B,\beta)$ is a $C$-morphism $f:A\rightarrow B$ such that $f\circ \alpha =\beta \circ F(f)$.
$$
\begin{tikzcd}[row sep=huge,column sep=huge]
F(A)
\arrow[r, "F(f)"]
\arrow[d, "\alpha"]
& F(B)
\arrow[d, "\beta"]
\\
A
\arrow[r, "f"]
& B
\end{tikzcd}
$$
\end{definition}

Given a signature $\Sigma$, we can define $F_\Sigma :\Set \rightarrow\Set$ as
$$
X \mapsto \sum_{f \in \Sigma} X^{|f|} \cong \sum_\N \Sigma_n \times X^n
$$
where $|f|$ denotes the arity of $f$, and identify each $\sigma$ with its corresponding summand $X^{|\sigma|}$, then the $\Sigma$-algebras are exactly the category of $F_\Sigma$-algebras, $\falg$. We can define a forgetful functor $U: \falg \rightarrow \Set$ by sending $(A, \alpha) \mapsto A$, which can be shown to have a left adjoint $T: V \mapsto T(V)$. The functor $T$ is fully faithful, because any homomorphism $T(V) \rightarrow T(V')$ is uniquely determined by its restriction to $V$ and the map associated to a homomorphism is its restriction to $V$. In particular the homomorphisms $T(V) \rightarrow T(V)$ are in bijection with maps $f: V \rightarrow T(V)$, where the homomorphism associated to $f$ is $T(f): T(V) \rightarrow T(f(V)) \subseteq T(V)$.

The initial algebras $T$ of $\falg$ are called algebraic data types (ADTs) in functional programming, and are used to give the following notion of a recursive program.

\begin{definition}[Catamorphism]
If $X$ is a set, $T$ is an initial $F$-algebra, and $c:F(X) \rightarrow X$, then the \textbf{catamorphism} of $c$ is the unique map $\cata(c) : T \rightarrow X$ such that
$$
\begin{tikzcd}[row sep=huge,column sep=huge]
F(T)
\arrow[r, "F(\cata(c))"]
\arrow[d, "\phi"]
& F(X)
\arrow[d, "c"]
\\
T
\arrow[r, "\cata(c)"]
& X
\end{tikzcd}
$$
commutes.
\end{definition}
 
For example, we can identify the natural numbers $\N$, along with the maps $* \mapsto 0$ and $n \mapsto n+1$, as the initial algebra of the endofunctor $X \mapsto X+1$. This lets us define recursive functions $\N \rightarrow X$ by specifying two functions $1 \rightarrow X$ and $X \rightarrow X$. For example the map $\alpha: \N \rightarrow \mathbb{C}$ which sends $n \mapsto z^n$ is given by the maps $* \mapsto 1$ and $w \mapsto zw$, since $\alpha(0) = \alpha \circ \phi (*) = 1$ and $\alpha(n+1) = \alpha \circ \phi (n) = z \alpha(n)$, and so $\alpha(n) = z^n \alpha(0) = z^n$ by induction.

We can dually define $F$-coalgebras for any endofunctor $F$, where instead we have morphisms $\alpha: A \rightarrow F(A)$. If we let $\psi: T \rightarrow F(T)$ by $f(t_1,...,t_n) \mapsto (t_1,...,t_n)$ then by Lambek's Lemma \cite{Lambek1968} $(T, \psi)$ is the final coalgebra, and $\psi$ is inverse to $\phi$.

We can push the algebra and coalgebra maps $\phi$ and $\psi$ through the functor $F$ to make $T$ into an $F^n$ algebra and coalgebra, and define
$$
\Phi_n := \phi \circ \ddd \circ F^{n-1}(\phi): F^n(T) \rightarrow T
$$
$$
\Psi_n := F^{n-1}(\psi) \circ \ddd \circ \psi: T \rightarrow F^n(T)
$$
which are mutual inverses. These both clearly satisfy
$$
\begin{tikzcd}[row sep=huge,column sep=huge]
F^n(T)
\arrow[r, "F^n(f)"]
\arrow[d, "\Phi_n"]
& F^n(T)
\arrow[d, "\Phi_n"]
\\
T
\arrow[r, "f"]
& T
\end{tikzcd}
\qquad
\begin{tikzcd}[row sep=huge,column sep=huge]
F^n(T)
\arrow[r, "F^n(f)"]
& F^n(T)
\\
T
\arrow[r, "f"]
\arrow[u, "\Psi_n"]
& T
\arrow[u, "\Psi_n"]
\end{tikzcd}
$$
for any homomorphism $f: T \rightarrow T$. The map $\Psi_n$ can be seen as splitting a (ground) term's expression tree into its constituent subtrees of distance at most $n$ from the root. Correspondingly, $\Phi_n$ will assemble a tree from a collection of subtrees. While the ground terms form an initial algebra and hence are also a coalgebra, if we take a set of terms with variables $T(V)$ then, while this is an $F$-algebra, it fails to be a coalgebra because the inverse map $\psi$ is not defined on the variables. We can, however, define it on the subset of non-variable terms, and by extension define $\Psi_n$ on all terms whose variables occur at depth at least $n$ from the root.

For example, if $x$, $y$, and $z$ are variables, the term $f(x,g(y,z))$ can be split by $\psi$ into $(x,g(y,z)) \in T^2 \subset F(T)$. However $\Psi_1 = F(\psi) \circ \psi$ cannot be applied since $x$ is a variable and so $\psi(x)$ is not defined.

So in the case that a term $t$ has all its leaves at the same depth, the right choice of $n$ will strip the tree down to just the variables and constants at the leaves, i.e. there exists some $n$ such that $\hat{t} := \Psi_n(t) \in (\Sigma_0 \cup V)^k \subset F^n(T)$ for some $k$. This will be important for our model as we would like to keep track of the action of a homomorphism on the tuple of variables in a given term. We can make this explicit and show that the correct choices of $n$ and $k$ are the \emph{depth} and \emph{leaf number} of the term, which we now define.

\begin{definition}
Let $t \in T(V)$ be a term. Then the \textbf{depth} of $t$ is defined as
$$
d(t) := \max \{length(p) : p \in Pos(t)\}
$$
and the \textbf{leaf number} of $t$ is defined as
$$
L(t) := |\{p \in Pos(t): t|_p \in V \cup \Sigma_0\}|.
$$
\end{definition}

Even when $t$ has no constants, the leaf number of is not necessarily the size of $Var(t)$, since variables in $t$ may be repeated more than once. For example if $t = f(g(x,y),x)$, where $x$ and $y$ are variables, then $d(t) = 2$ and $L(t) = 3$.

\subsection{Reduction Relations and Rewriting Functions}

We can now make use of the variables to define expressions in the term algebra. These will be used to parameterise rewriting processes, which can be induced by \emph{identities} of two terms.

\begin{definition}[$\Sigma$-identity]
Let $\Sigma$ be a signature and $V$ a set of variables disjoint from $\Sigma$. A $\Sigma$\textbf{-identity} (or simply identity) is a pair $(s,t) \in T(V) \times T(V)$. We call $s$
the left-hand side (LHS) and $t$ the right-hand side (RHS), and assume always that $Var(t) \subseteq Var(s)$.
\end{definition}

Identities are the basis for rewriting, by generating from each identity a \emph{reduction relation}. This nomenclature refers to the typical use of rewriting as a simplification of complex expressions. In our application, the relation will in fact be building up a sequence of increasingly complex terms, and so is better viewed as an expansion rather than a reduction.

\begin{definition}[Reduction Relation]
Let E be a set of $\Sigma$-identities. The reduction relation
$\rightarrow_E \subset T(V) \times T(V)$ is defined as
$s \rightarrow_E t$ if there exists some pair $(l,r) \in E$, $p \in Pos(s)$, and $\sigma \in Sub$ where
$$
s|_p = \sigma(l) \textnormal{ and } t = s[\sigma(r)]_p.
$$
We sometimes write $s \rightarrow_E^p t$ to indicate the position at which the reduction takes place. We also call $(l,r)$ a \textbf{rewrite rule} and say that $s$ \textbf{rewrites} to $t$.
\end{definition}

For example, if
$$
l = f(x,f(y,z))
\qquad
r = f(f(x,y),z)
$$
is an identity, then $f(a,f(g(b),f(c,d)))$ will rewrite to $f(f(a,g(b)),f(c,d))$ at position $\epsilon$ and rewrite to $f(a,f(f(g(b),c),d))$ at position 2.

Reduction relations are typically used in rewriting systems comprising multiple rules and where the rewriting can occur at any position. For our application, we will be using a single rewrite rule (although in practice this may be searched for and constructed as a composition of multiple simple rules), and insist that the reduction occurs at a specified position. In this case, the reduction relation defines a function between the subsets of $T(V)$ which \emph{match} with the left-hand rule. We note for our application that the set of ground terms $T$ is clearly closed under rewriting, since $Var(t) \subseteq Var(s)$ whenever $s \rightarrow_E t$.

\begin{definition}
Given a term $s$ and position $p$, we define the set of \textbf{instances of $s$ at position $p$} by
$$
T_p^s := \{t \in T(V): t|_p = \sigma(s) \textnormal{ for some } \sigma \in Sub\}.
$$
Given a pair of terms $(l,r)$, where $Var(r) \subseteq Var(l)$, we can define a \textbf{rewriting function} $R_p: T_p^l \rightarrow T_p^r$ by
$$
t \mapsto t[\sigma(r)]_p
$$
if $t|_p = \sigma(l)$.
\end{definition}

We need to check that $R_p$ is well defined, in the sense that it does not depend on the substitution $\sigma$. We further show that it is a surjection, and is also injective in the case that $Var(r) = Var(l)$.

\begin{lemma}
$R_p$ is a well defined surjection $T_p^l \rightarrow T_p^r$ which is injective if $Var(r) = Var(l)$.
\end{lemma}
\begin{proof}
In order to show that $R_p$ is well defined, we need to check that if $\sigma(l) = \sigma'(l)$ for two substitutions $\sigma$ and $\sigma'$, then $\sigma(r) = \sigma'(r)$. We claim that if $\sigma(t) = \sigma'(t)$ for some term $t$ then $\sigma = \sigma'$ on $Var(t)$. If $v \in Var(t)$ then there is some $p$ where $t|_p = v$. So
$$
\sigma(v) = \sigma(t|_p) = \sigma(t)|_p = \sigma'(t)|_p = \sigma'(t|_p) = \sigma'(v).
$$
We also note by induction that $\sigma(t)$ is determined by the value of $\sigma$ on $Var(t)$. Using this and the fact that $Var(r) \subseteq Var(l)$, it follows that $\sigma(r) = \sigma'(r)$ whenever $\sigma(l) = \sigma'(l)$.

It is clear that the codomain of $R_p$ is indeed contained in $T_p^r$, since $(t[\sigma(r)]_p)|_p = \sigma(r)$ and so if $t \in T_p^l$ then $R_p(t) \in T_p^r$. Conversely, if $t \in T_p^r$ then $t|_p = \sigma(r)$ for some $\sigma$. Now let $t' = t[\sigma(l)]_p \in T_p^l$. We see that
$$
R_p(t') = \big(t[\sigma(l)]_p \big) [\sigma(r)]_p = t[\sigma(r)]_p = t,
$$
so $R_p$ is surjective.

We now show that $R_p$ is injective under the additional assumption that $Var(r) = Var(l)$. Suppose $t, t' \in T_p^l$, so $t|_p = \sigma(l)$ and $t'|_p = \sigma'(l)$ for substitutions $\sigma$ and $\sigma'$, and that $R_p(t) = R_p(t')$ i.e. $t[\sigma(r)]_p = t'[\sigma'(r)]_p$.

We can see that
$$
\sigma(r) = \big(t[\sigma(r)]_p \big)|_p = \big(t'[\sigma'(r)]_p \big)|_p = \sigma'(r).
$$
By the same argument as above, where this time we use the fact that $Var(r) = Var(l)$, we have $\sigma(l) = \sigma'(l)$ and so $t|_p = t'|_p$. It then follows that
\begin{equation*}
    \begin{split}
        t &= t[t|_p]_p \\
        &= \big(t[\sigma(r)]_p \big)[t|_p]_p \\
        &= \big(t'[\sigma'(r)]_p \big)[t|_p]_p \\
        &= t'[t|_p]_p \\
        &= t'[t'|_p]_p \\
        &= t',
    \end{split}
\end{equation*}
and so $R_p$ is injective.
\end{proof}

\begin{example}
Suppose $\Sigma$ is a signature comprising two constants $a$ and $b$, as well as a single binary relation $(\cdot, \cdot)$. Suppose we introduce two variables $x$ and $y$ to define a rewrite rule $(x,y) \mapsto (y, (x,y))$. Then, for example,
$$
T_\epsilon^l = \{(t_1, t_2): t_i \in T(V)\}
$$
and
$$
T_{21}^l = \{(t_1, ((t_2, t_3), t_4)): t_i \in T(V)\}.
$$
The map $R_{21}$ sends
$$
(t_1, ((t_2, t_3), t_4)) \mapsto (t_1, ((t_3, (t_2, t_3)), t_4)).
$$
\end{example}

The requirement that $Var(r) = Var(l)$ for injectivity is essential. If we take $l = f(x,y)$ and $r = g(x)$ for variables $x$ and $y$ then $R_\epsilon$ is clearly not injective, as $R_\epsilon(f(x,y)) = R_\epsilon(f(x,x)) = g(x)$.

For our application, we are interested in rewrite rules that can be repeatedly applied to a single term to generate a sequence. This is possible exactly when $r = \tau(l)$ for some substitution $\tau$.

\begin{lemma}\label{iterable_rr}
$T_p^r \subseteq T_p^l$ if and only if $r = \tau(l)$ for some substitution $\tau$.
\end{lemma}
\begin{proof}
If $r = \tau(l)$ for some $\tau$, and $t \in T_p^r$, then $t|_p = \sigma(r) = \sigma\tau(l)$ for some $\sigma$. So $t \in T_p^l$ and hence $T_p^r \subseteq T_p^l$. Conversely, suppose $T_p^r \subseteq T_p^l$. Let $t$ be any term with $t|_p = r$, so that $t \in T_p^r$ is witnessed by the identity substitution. Then by assumption $t \in T_p^l$, so $r = t|_p = \tau(l)$ for some $\tau$.
\end{proof}

We would like to establish a connection between an iterated rewriting function and the \emph{hidden} state of a recurrent system, which can be in any type. We do this by associating to each reduction relation a natural transformation between powers of the endofunctor $F$, such that the rewriting function $R_p$ is conjugate to the action of this natural transformation on $F(T)$.

$$
\begin{tikzcd}[row sep=huge,column sep=huge]
F^{d(l)}(T(V))
\arrow[r, "\eta_{T(V)}"]
& F^{d(r)}(T(V))
\arrow[d, "\Phi_{d(r)}"]
\\
T_l^p
\arrow[r, "R_\epsilon"]
\arrow[u, "\Psi_{d(l)}"]
& T_l^p
\end{tikzcd}
$$

If we consider a dynamical system on $\R^2$ as an example, where the hidden dynamical process is given by $g = (g_1,g_2)$ and the output function by $f$, we can represent the symbolic application of $g$ by the rewrite rule $f(x_1,x_2) \mapsto f(g_1(x_1,x_2), g_2(x_1,x_2))$, where $x_1$ and $x_2$ are variables. We can then find lifted terms $\hat{l} \in V^2 \subset F(T)$ and $\hat{r} \in V^4 \subset F^2(T)$. Viewing $F$ as a polynomial functor in the variable $Y$, we can use this information to define a natural transformation between the corresponding cofactors $Y^2$ and $Y^4$ of $F$ and $F^2$, by defining the map $4 \rightarrow 2$ as $(1,2,1,2)$ determined by the variables $x_1$ and $x_2$. If we can somehow extend this natural transformation to the other cofactors of $F$, we can define $\eta: F \rightarrow F^2$ such that $\hat{r} = \eta_T(\hat{l})$. The naturality of $\eta$ will allow us to define an analogous \emph{rewrite} on other $F$-algebras. In this example the map $\eta_\R: \R^2 \rightarrow \R^4$ will describe the application of $g$ and $f$.

There is, however, a technical problem with generalising the case above to all examples of rewrite rules. For example if the rule was $f(x_1,x_2) \mapsto f(g_1(x_1,x_2), x_2)$ then, although $\hat{r} \in V^3 \subset F^2(T)$, there is no copy of $Y^3$ in $F^2$ with which to identify this, because the subterms of $r$ have different depths. We can, however, account for this by including an additional unary operator $\iota$ in $\Sigma$ which will always be interpreted as the identity $Y \rightarrow Y$ for any $\Sigma$-algebra. In this example, we can then write $r = f(g_1(x_1,x_2), \iota(x_2))$, so that the rewrite corresponds to a natural transformation of cofactors $Y^2 \rightarrow Y^3$. We make this precise below, where we write $\Sigma'$ for this extended signature, but first define the property of all variables occurring at the same depth, which we call \emph{flat}.

\begin{definition}
A term $t \in T(V)$ is called \textbf{flat} if $l(p) = d(t)$ for all $p \in Pos(t)$ with $t|_p \in V$.
\end{definition}

We now check that the extension of $\Sigma$ to $\Sigma'$ guarantees for each term $t$ the existence of an essentially equivalent flat term $t'$. If $(Y,\alpha)$ is a $\Sigma$-algebra we let $(Y,\alpha')$ be the $\Sigma'$-algebra where $\alpha'$ extends $\alpha$ by $\alpha': \iota \mapsto id_Y$

\begin{lemma}\label{equal_depth_lemma}
If $t \in T(V)$, then there is a flat term $t' \in T(\Sigma',V)$ such that $\cata(\alpha)(t) = \cata(\alpha')(t')$ for all $\Sigma$-algebras $(Y,\alpha)$.
\end{lemma}
\begin{proof}
We define $t'$ by induction on $d(t)$. If $d(t) = 0$ then $Pos(t) = \{\epsilon\}$ so $t' = t$. If $d(t) \geq 1$ then $t = f(t_1,...,t_n)$ for some $f \in \Sigma_n$ (and where $d(t_i)\leq d(t)-1$ for all $i$) and so we set
$$
t' = f(\iota^{d(t)-d(t_1)-1}(t_1'),...,\iota^{d(t)-d(t_n)-1}(t_n')).
$$
Then, if $t'|_p \in V$, we have $p = iqp'$ where $1 \leq i \leq n$, $q$ is a string of $d(t)-d(t_i)-1$ ones, and $q \in Pos(t_i')$ where $t_i'|_q = t'|_p \in V$. So $l(q) = d(t_i)$ by induction, and $l(p) = 1 + (d(t)-d(t_i)-1) + d(t_i) = d(t)$, so $t'$ is flat.

We can then check inductively that
\begin{equation*}
    \begin{split}
        \cata(\alpha')(t')
        &= \cata(\alpha')(f(\iota^{d(t)-d(t_1)-1}(t_1'),...,\iota^{d(t)-d(t_n)-1}(t_n'))) \\
        &= \alpha(f)(\cata(\alpha')(\iota^{d(t)-d(t_1)-1}(t_1')),...,\cata(\alpha')(\iota^{d(t)-d(t_n)-1}(t_n'))) \\
        &= \alpha(f)(\cata(\alpha')(t_1'),...,\cata(\alpha')(t_n')) \qquad \textnormal{since } \alpha'(\iota) = id_Y \\
        &= \alpha(f)(\cata(\alpha)(t_1),...,\cata(\alpha)(t_n)) \qquad \textnormal{by induction} \\
        &=\cata(\alpha)(f(t_1,...,t_n)) \\
        &=\cata(\alpha)(t) \\
    \end{split}
\end{equation*}
for all $\Sigma$-algebras $(Y,\alpha)$.
\end{proof}

So with the addition of this additional identity operator, we can assume without loss of generality that the left and right-hand side terms of the rewrite rule are flat. This is precisely the condition that guarantees we can identify the variable tuple $V^{L(t)}$ with a cofactor $Y^{L(t)}$ of $F^{d(t)}$. As noted before we can define the map $\Psi_n$ on all terms with minimal leaf depth at least $n$, so in particular we can lift each flat term $t$ to a term $\hat{t} := \Psi_n(t)$. We now make this more precise.

\begin{lemma}\label{lifted_terms_lemma}
Let $t \in T(V)$ be flat. Then the functor $Y \mapsto F^{d(t)}(Y)$ has a cofactor $Y^{L(t)}$, and there is a tuple $\hat{t} \in (\Sigma_0 \cup V)^{L(t)}$ such that $\Phi_{d(t)}(\hat{t}) = t$.
\end{lemma}
\begin{proof}
We use induction on $d(t)$, and use $Y^n$ to denote the monomial functor throughout. If $d(t) = 0$, then $t$ is either a constant or a variable, meaning $L(t) = 1$. It follows immediately that $F^0 = Y$ has $Y^1 = Y$ as a cofactor, and that $t \in \Sigma_0 \cup V \subset T(V)$ satisfies $\Phi_0(t) = t$.

If $d(t) \geq 1$ then $t = f(t_1,...,t_n)$ for some $f \in \Sigma_n$, where all $t_i$ have $d(t_i) = d(t) - 1$ (since flatness is a hereditary property). Note that we can condition
$$
Pos(t) = \bigcup_{i=1}^n \{ip' : p' \in Pos(t_i)\}
$$
and so
$$
\{p \in Pos(t): t|_p \in \Sigma_0 \cup V\} = \bigcup_{i=1}^n \{ip': p' \in Pos(t_i), t|_p \in \Sigma_0 \cup V\}
$$
is a disjoint union, giving
$$
L(t) = \sum_{i=1}^n |\{p' \in Pos(t_i): t|_p \in \Sigma_0 \cup V\}| = \sum_{i=1}^n L(t_i).
$$
If we make the inductive assumption that $Y \mapsto F^{d(t)-1}(Y)$ has $Y^{L(t_i)}$ as a cofactor for all $i$ then $F^{d(t)}(Y) = F(F^{d(t)-1}(Y))$ has
$$
\prod_{i=1}^n Y^{L(t_i)} \cong Y^{\sum_{i=1}^n L(t_i)} \cong Y^{L(t)}
$$
as a cofactor. We can then define
$$
\hat{t} = (\hat{t_1},...,\hat{t_n}) \in \prod_{i=1}^n (\Sigma_0 \cup V)^{L(t_i)} \cong (\Sigma_0 \cup V)^{L(t)} \subset F^{d(t)}(\Sigma_0 \cup V),
$$
and verify that
\begin{equation*}
    \begin{split}
        \Phi_{d(t)}(\hat{t})
        &= \phi \circ F(\Phi_{d(t)-1})(\hat{t}) \\
        &= \phi (\Phi_{d(t)-1}(\hat{t_1}),...,\Phi_{d(t)-1}(\hat{t_n})) \\
        &= \phi (t_1,...,t_n) \\
        &= t
    \end{split}
\end{equation*}
by induction.
\end{proof}

If we take the example $t = f(g(f(x,c),y),x)$, where $c$ is a constant and $x,y$ are variables, then $d(t) = 3$ and $L(t) = 4$. We can flatten $t$ to the $T(\Sigma', V)$ term $t' = f(g(f(x,c),\iota(y)),\iota^2(x))$, so that $\hat{t'} = (x,c,y,x) \in (\Sigma_0 \cup V)^4$. We can check that
\begin{equation*}
    \begin{split}
        \Phi_3(\hat{t'})
        &= f(\Phi_2(x,c,y),\Phi_2(x) \\
        &= f(g(\phi(x,c),\phi(y)),\iota(\phi(x))) \\
        &= f(g(f(x,c),\iota(y)),\iota^2(x)) \\
        &= t'
    \end{split}
\end{equation*}
according to the configuration of the cofactor $Y^4$ in $F^3$. The representable notation $Y^n$ does obscure the fact that when we refer to such terms as powers of $F^d$ they also contain the information about the construction of the terms within. For example, there is a copy of $Y^2$ in $F^2$ corresponding to both terms of the form $g(f(x,y))$ and $f(g(x),y)$, although this sort of algebraic ambiguity will not be important in the following application.

Before proceeding with the result we summarise some basic facts about polynomial functors. For a full exposition see \cite{2005.01894}.

\begin{lemma}[Yoneda]
Given a functor $F : \Set \rightarrow \Set$ and a set $S$, there is an isomorphism
$$
F(S) \cong Nat(Y^S, F)
$$
where $Nat$ denotes the set of natural transformations. Moreover, this isomorphism is natural in
both $S$ and $F$.
\end{lemma}

In particular, we obtain the useful corollary that
$$
S^T \cong Nat(Y^S, Y^T)
$$
which classifies the monomial natural transformations. In the case of polynomials, we use the fact that the natural transformations on coproduct functors are exactly the coproducts of natural transformations on each cofactor, so that, if $p$ and $q$ are polynomials,
$$
Nat(p, q) = \prod_{i\in p(1)} \sum_{j\in q(1)} p_i^{q_j}.
$$

We can now prove a lemma that associates to each rewrite rule a natural transformation that factors through the algebra structure. Intuitively we would like to take a term $t$ that matches with a flat left-hand term $l$ at position $\epsilon$, and \emph{push} $t$ up the levels of $F$ to depth $d(l)$, corresponding to breaking the expression tree for $t$ down into its subtrees at the variable positions of $l$. The rewriting of these variables from $l$ to $r$ can then be conferred onto $t$ as a natural transformation from $T(V)^{Var(l)} \rightarrow T(V)^{Var(r)}$. In practice, it will prove more straightforward to extend this to a natural transformation $\eta: T(V)^{d(l)} \rightarrow T(V)^{d(r)}$. Conjugation with $\Psi_{d(l)}$ and $\Phi_{d(r)}$ will allow us to express the rewriting function in terms of $\eta_{T(V)}$. We will later deal with the case of non-trivial rewrite position by restricting the term to the subterm at position $p$ and applying a position-$\epsilon$ rewrite, and so this simple case is all we characterise here.

Note that we cannot deal with non-trivial rewrite position by simply pushing this same process up to the depth of the rewrite position. This is because the term may match multiple times at the same depth, but we only want the subterm at the specified position to change. In other words $\Phi_{l(p)} \circ R_\epsilon \circ \Psi_{l(p)} \neq R_p$.

\begin{lemma}\label{nat_trans_thm}
If $(l,r)$ is a $\Sigma'$-identity, where $l$ and $r$ are flat and $d(r) \geq d(l)$, then there exists a natural transformation $\eta: F^{d(l)} \rightarrow F^{d(r)}$ such that the diagram
$$
\begin{tikzcd}[row sep=huge,column sep=huge]
F^{d(l)}(T(V))
\arrow[r, "\eta_{T(V)}"]
& F^{d(r)}(T(V))
\arrow[d, "\Phi_{d(r)}"]
\\
T_l^p
\arrow[r, "R_\epsilon"]
\arrow[u, "\Psi_{d(l)}"]
& T_l^p
\end{tikzcd}
$$
commutes for all positions $p$.
\end{lemma}
\begin{proof}
Given the two flat terms $l, r \in T(V)$, by Lemma \ref{lifted_terms_lemma} we can find tuples
$$
\hat{l} = \Psi_{d(l)}(l) \in (\Sigma_0 \cup V)^{L(l)} \subset F^{d(l)}(T(V))
$$
$$
\hat{r} = \Psi_{d(r)}(r) \in (\Sigma_0 \cup V)^{L(r)} \subset F^{d(r)}(T(V)).
$$
We begin by constructing $\eta: F^{d(l)} \rightarrow F^{d(r)}$ such that $\hat{r} = \eta_{T(V)} (\hat{l})$. Since we assume that $Var(r) \subseteq Var(l)$, we can assign to each $i \in \{1,...,L(r)\}$ some $j_i \in \{1,...,L(l)\}$ such that $\hat{r}_i = \hat{l}_{j_i}$. So let $u:L(r) \rightarrow L(l)$ by $i \mapsto j_i$, which defines a natural transformation $Y^{L(l)} \rightarrow Y^{L(r)}$. To extend this to a natural transformation $\eta: F^{d(l)} \rightarrow F^{d(r)}$ we must pick natural transformations on the other cofactors of $F^{d(l)}$. Under the assumption that $d(r) \geq d(l)$ we can use the $d(r) - d(l)$ power of the identity operator $\iota$ to match each $Y^n\subset F^{d(l)}(Y)$ with $Y^n\subset F^{d(r)}(Y)$. We choose the identity natural transformation so that $\eta: F^{d(l)} \rightarrow F^{d(r)}$ will preserve the cofactors that are not rewritten. By construction we have that $\hat{r} = \eta_{T(V)} (\hat{l})$.

From the properties of $\Phi_n$ and $\Psi_n$, and naturality of $\eta$, we can assemble a commutative diagram
$$
\begin{tikzcd}[row sep=huge,column sep=huge]
T(V)
\arrow[r, "\Psi_{d(l) }"]
\arrow[d, "\sigma"]
& F^{d(l)}T(V)
\arrow[r, "\eta_{T(V)}"]
\arrow[d, "F^{d(l)}(\sigma)"]
& F^{d(r)}T(V)
\arrow[r, "\Phi_{d(r)}"]
\arrow[d, "F^{d(r)}(\sigma)"]
& T(V)
\arrow[d, "\sigma"]
\\
T(V)
\arrow[r, "\Psi_{d(l)}"]
& F^{d(l)}T(V)
\arrow[r, "\eta_{T(V)}"]
& F^{d(r)}T(V)
\arrow[r, "\Phi_{d(r)}"]
& T(V)
\end{tikzcd}
$$

If $s \in T_\epsilon^l$ and $t = R_\epsilon (s)$, then there exists some substitution $\sigma$ with $s = \sigma(l)$ and $t = s[\sigma(r)]_\epsilon = \sigma(r)$. We can evaluate the diagram at $l$ to obtain
$$
\begin{tikzcd}[row sep=huge,column sep=huge]
l
\arrow[r, "\Psi_{d(l)}"]
\arrow[d, "\sigma"]
& \Psi_{l(p)}(\hat{l})
\arrow[r, "\eta_{T(V)}"]
& \Psi_{l(p)}(\hat{r})
\arrow[r, "\Phi_{d(r)}"]
& r
\arrow[d, "\sigma"]
\\
s
\arrow[rrr, "\Phi_{d(r)}\circ \eta_{T(V)}\circ \Psi_{d(l)}", dashed]
&&& t
\end{tikzcd}
$$
and so $t = \Phi_{d(r)}\circ \eta_{T(V)}\circ \Psi_{d(l)}(s)$. Conversely suppose that $s \in T_\epsilon^l$ (so $s = \sigma(l)$ for some substitution $\sigma$) and $t = \Phi_{d(r)}\circ \eta_{T(V)}\circ \Psi_{d(l)}(s)$. We now have
$$
\begin{tikzcd}[row sep=huge,column sep=huge]
l
\arrow[r, "\Psi_{d(l)}"]
\arrow[d, "\sigma"]
& \Psi_{l(p)}(\hat{l})
\arrow[r, "\eta_{T(V)}"]
& \Psi_{l(p)}(\hat{r})
\arrow[r, "\Phi_{d(r)}"]
& r
\arrow[d, "\sigma", dashed]
\\
s
\arrow[rrr, "\Phi_{d(r)}\circ \eta_{T(V)}\circ \Psi_{d(l)}"]
&&& t
\end{tikzcd}
$$
so $t = \sigma(r) = R_\epsilon(s)$. So $R_\epsilon = \Phi_{d(r)}\circ \eta_{T(V)}\circ \Psi_{d(l)}$ for all positions $p$.
\end{proof}

Note that the natural transformation $\eta$ is not uniquely determined by the pair $(l,r)$, since the term $l$ may have repeated variables. In the case that each variable in $l$ appears only once, and $l$ contains no constants, then there is a unique $\eta$.

\section{Rewriting Models}

We will now introduce rewriting models for time series, and prove that they are indeed the algebraic analogue of dynamical systems.

We start by defining a term algebra on a signature, and some rewrite rule in that algebra. We would like the rewrite rule to represent some time-homogeneous process, so we will generate a sequence from the iterated application of that single rule to some initial term $t_0$. By Lemma \ref{iterable_rr} this is possible precisely when $r = \tau(l)$. We also choose a $\Sigma$-algebra $(X,c)$ which comprises the output set $X$ and the map $c$ which interprets the terms as elements of $X$. From these constituent pieces, we can define our model.

\begin{definition}
Let $\Sigma$ be a signature, $(l,r')$ be a $\Sigma$-identity in some $T(V)$ with $r = \tau(l)$, and $p$ a position, which induces a rewriting function $R_p: T_p^l \rightarrow T_p^l$. If $(X,c)$ is a $\Sigma$-algebra, and $t_0 \in T_p^l$ is a ground initial term, then a \textbf{rewriting model} is the map $\N \rightarrow X$ given by $n \mapsto \cata(c) \circ R_p^n(t_0)$.
\end{definition}

Note that, if we require that $r = \tau(l)$, then if $l$ is flat then $r$ is not flat in general. This is because substitution increases the leaf depth of the variables but not the constants. So even if $\tau$ only substitutes flat terms of the same depth, $\tau(l)$ will still not be flat if $l$ contains constants. Hence we define $r$ as a general term but rewrite using the flattened term $r'$.

We will show that this model is equivalent to a dynamical system on $X^{L(l)}$. First, we define two maps
$$
H: F^{d(l)}(T) \rightarrow F^{d(l)}(T) \qquad H = F^{d(r)}(\phi) \circ \ddd \circ F^{{d(l)}-1}(\phi) \circ \eta_T
$$
and
$$
G: F^{d(l)}(X) \rightarrow F^{d(l)}(X) \qquad \Phi = F^{d(l)}(c) \circ \ddd \circ F^{{d(r)}-1}(c) \circ \eta_X,
$$
where $\eta$ is the induced natural transformation $F^{d(l)} \rightarrow F^{d(r)}$ in Theorem \ref{nat_trans_thm}. Note that the requirement that $d(r) \geq d(l)$ is satisfied by the fact that $r = \tau(l)$. $G$ will be the hidden dynamical process on $X^{L(l)}$, and $H$ is the algebraic equivalent on $T^{L(l)}$.

\begin{lemma}\label{cat_lemma}
Let $H$ and $G$ be as above. Then
$$
F^{d(l)}(\cata(c)) \circ H^n = G^n \circ F^{d(l)}(\cata(c)),
$$
and this also holds on the cofactors $T^{L(l)}$ and $X^{L(l)}$.
\end{lemma}
\begin{proof}
It suffices by induction to prove the case $n=1$. Naturality of $\eta: F^{d(l)} \rightarrow F^{d(r)}$ gives the commutative square
$$
\begin{tikzcd}[row sep=huge,column sep=huge]
F^{d(l)}(T)
\arrow[r, "\eta_T"]
\arrow[d, "F^{d(l)}(\cata(c))"]
& F^{d(r)}(T)
\arrow[d, "F^{d(r)}(\cata(c))"]
\\
F^{d(l)}(X)
\arrow[r, "\eta_X"]
& F^{d(r)}(X)
\end{tikzcd}
$$
Inductively applying $F$ to the definition of catamorphism produces the ladder
$$
\begin{tikzcd}[row sep=huge,column sep=normal]
F^{d(r)}(T)
\arrow[r, "F^{d(r)-1}(\phi)"]
\arrow[d, "F^{d(r)}(\cata(c))"]
&
F^{d(r)-1}(T)
\arrow[r, "F^{d(r)-2}(\phi)"]
\arrow[d, "F^{d(r)-1}(\cata(c))"]
&
\ddd
\arrow[r, "F^{d(l)+1}(\phi)"]
&
F^{d(l)+1}(T)
\arrow[r, "F^{d(l)}(\phi)"]
\arrow[d, "F^{d(l)+1}(\cata(c))", swap]
&
F^{d(l)}(T)
\arrow[d, "F^{d(l)}(\cata(c))", swap]
\\
F^{d(r)}(X)
\arrow[r, "F^{d(r)-1}(c)"]
&
F^{m-1}(X)
\arrow[r, "F^{d(r)-2}(c)"]
&
\ddd
\arrow[r, "F^{d(l)+1}(c)"]
&
F^{d(l)+1}(X)
\arrow[r, "F^{d(l)}(c)"]
&
F^{d(l)}(X)
\end{tikzcd}
$$
where $d(r) \geq d(l)$ because $r = \tau(l)$. Gluing these along $F^{d(r)}(\cata(c))$ gives the result. 
$$
\begin{tikzcd}[row sep=huge,column sep=huge]
F^{d(l)}(T)
\arrow[r, "\eta_T"]
\arrow[d, "F^{d(l)}(\cata(c))"]
\arrow[rrr, "H", bend left=20]
& F^{d(r)}(T)
\arrow[rr, "F^{d(l)}(\phi) \circ \ddd \circ F^{d(r)-1}(\phi)"]
\arrow[d, "F^{d(r)}(\cata(c))"]
&& F^{d(l)}(T)
\arrow[d, "F^{d(l)}(\cata(c))"]
\\
F^{d(l)}(X)
\arrow[r, "\eta_X"]
\arrow[rrr, "G", bend right=20]
& F^{d(r)}(X)
\arrow[rr, "F^{d(l)}(c) \circ \ddd \circ F^{d(r)-1}(c)"]
&& F^{d(l)}(X)
\end{tikzcd}
$$
We constructed $\eta$ such that $\eta: Y^{L(l)} \rightarrow Y^{L(r)}$ (where $Y^n$ denotes the functor $Y \mapsto Y^n$), so we obtain the following restriction.
$$
\begin{tikzcd}[row sep=huge,column sep=huge]
T^{L(l)}
\arrow[r, "\eta_T"]
\arrow[d, "\cata(c)^{L(l)}"]
& T^{L(r)}
\arrow[d, "\cata(c)^{L(l)}"]
\\
X^{L(l)}
\arrow[r, "\eta_X"]
& X^{L(r)}
\end{tikzcd}
$$
We now check if $r = \tau(l)$ then the ladder construction is also valid. Since $l$ is flat it lifts to the tuple $\hat{l} \in (V \cup \Sigma_0)^{L(l)} \subset F^{d(l)}(T)$. If $r = \tau(l)$ then $\Psi_{d(l)}(r) \in T^{L(l)}$ also, but each variable in $\hat{l}$ is now replaced with its substitution under $\tau$. Flattening to $r'$ means that each component of the tuple $\Psi_{d(l)}(r') \in T^{L(l)}$ has the same depth. Note that, although it is suppressed by the notation $Y^n$, each representable cofactor of $F^k$ contains the compositional information about how they were constructed from lower-order terms. As such we obtain
$$
\begin{tikzcd}[row sep=huge,column sep=huge]
T^{d(r)}
\arrow[rr, "F^{d(l)}(\phi) \circ \ddd \circ F^{d(r)-1}(\phi)"]
\arrow[d, "\cata(c)^{d(r)}"]
&& T^{d(l)}
\arrow[d, "\cata(c)^{d(l)}"]
\\
X^{d(r)}
\arrow[rr, "F^{d(l)}(c) \circ \ddd \circ F^{d(r)-1}(c)"]
&& X^{d(l)}
\end{tikzcd}
$$
since all maps involved respect the algebra structure. So we have
$$
\begin{tikzcd}[row sep=huge,column sep=huge]
T^{L(l)}
\arrow[d, "\cata(c)^{d(l)}"]
\arrow[r, "H"]
& T^{L(l)}
\arrow[d, "\cata(c)^{d(l)}"]
\\
X^{L(l)}
\arrow[r, "G"]
& X^{L(l)}
\end{tikzcd}
$$
where $G$ and $H$ denote the relevant restrictions.
\end{proof}

This lemma tells us that the inductive definition of the term sequence has a corresponding dynamical process in $X^{d(l)}$, and these are conjugate with the lifted catamorphism $F^{d(l)}(\cata(c))$. This is at the core of the following main result, which says that the rewriting models are an algebraically enriched class of dynamical systems. Specifically, the cartesian dynamical systems embed into the rewriting models, and each rewriting model projects onto a dynamical system, where that projection admits a section. So we can think of the rewriting models as a space of models that is fibred over the dynamical systems, where the fibre contains the additional structural information.

\begin{theorem}[Equivalence of models]\label{main_thm}
The class of cartesian dynamical systems embeds in the class of rewriting models. The rewriting models project onto the cartesian dynamical systems and that projection admits a section.
\end{theorem}
\begin{proof}
We will show how, given a model in one of these classes, we can construct a corresponding model in the other. We will see that the composition of these constructions will fix any dynamical system, and so give an embedding. Conversely, the composed constructions may forget the additional algebraic structure given by a rewriting model.  We will first prove that any rewriting model will project onto a cartesian dynamical system, and then show that this projection has a section (or left-inverse).

We start by considering a rewriting model at the base position $\epsilon$, given by $(l,r)$ and $c$, and some initial term $t_0$. As in Lemma \ref{cat_lemma} we can construct a ladder using the definition of catamorphism.
$$
\begin{tikzcd}[row sep=huge,column sep=large]
F^{v(l)}(T)
\arrow[r, "F^{v(l)-1}(\phi)"]
\arrow[d, "F^{v(l)}(\cata(c))"]
&
F^{v(l)-1}(T)
\arrow[r, "F^{v(l)-2}(\phi)"]
\arrow[d, "F^{v(l)-1}(\cata(c))"]
&
\ddd
\arrow[r, "F(\phi)"]
&
F(T)
\arrow[r, "\phi"]
\arrow[d, "F(\cata(c))"]
&
T
\arrow[d, "\cata(c)"]
\\
F^v(l)(X)
\arrow[r, "F^{v(l)-1}(c)"]
&
F^{v(l)-1}(X)
\arrow[r, "F^{v(l)-2}(c)"]
&
\ddd
\arrow[r, "F(c)"]
&
F(X)
\arrow[r, "c"]
&
X
\end{tikzcd}
$$
Using Lemma \ref{cat_lemma} we can form the diagram
$$
\begin{tikzcd}[row sep=huge,column sep=huge]
T_\epsilon^l
\arrow[r, "\Psi_{v(l)}"]
\arrow[rrr, "R_\epsilon^n", bend left=25]
& T^{v(l)}
\arrow[r, "H^n"]
\arrow[d, "F^{v(l)}(\cata(c))"]
& T^{v(l)}
\arrow[d, "F^{v(l)}(\cata(c))"]
\arrow[r, "\Phi_{v(l)}"]
& T_\epsilon^l
\arrow[d, "\cata(c)"]
\\
& X^{v(l)}
\arrow[r, "G^n"]
& X^{v(l)}
\arrow[r, "C_{v(l)}"]
& X
\end{tikzcd}
$$
where $C_k := c \circ \ddd \circ F^{k-1}(c)$. If we define an initial state
$$
x_0 = F^{v(l)}(\cata(c))(\Psi_{v(l)}(t_0)) \in X^{v(l)}
$$
then
$$
\cata(c) \circ R_\epsilon^n(t_0) = C_{v(l)} \circ G^n(x_0),
$$
which is a cartesian dynamical system with internal state $X^{v(l)}$.

We now deal with the case in which the rewriting occurs at a non-trivial position $p$. To do this we would like to restrict our terms to their subterms at position $p$, which reduces the situation to the above case for position $\epsilon$. The technicality is that the output of the induced dynamical system will be the value of the rewritten subterm, and not the entire term we started with. We will need to define a \emph{context} function $X \rightarrow X$ which compensates for this adjustment. We first check that the restriction we want makes sense, so let $r_p: T_p^l \rightarrow T_\epsilon^l$ by $t \mapsto t|_p$. We would like the following diagram to commute
$$
\begin{tikzcd}[row sep=huge,column sep=huge]
T_p^l
\arrow[r, "R_p"]
\arrow[d, "r_p"]
& T_p^r
\arrow[d, "r_p"]
\\
T_\epsilon^l
\arrow[r, "R_\epsilon"]
& T_\epsilon^r
\end{tikzcd}
$$
and can verify that, if $t \in T_p^l$ so $t|_p = \sigma(l)$, we have
$$
r_p \circ R_p (t) = r_p(t[\sigma(r)]_p) = \sigma(r) = R_\epsilon(\sigma(l)) = R_\epsilon \circ r_p(t).
$$
We now construct our function $X \rightarrow X$ which will have the effect of substituting the output of the dynamical system back into the initial term. We define a substitution function $Sub_s^p : T(V) \rightarrow T(V)$ by
$$
Sub_s^p(t) =
\left\{
	\begin{array}{ll}
		s[t]_p  & \mbox{if } p \in Pos(s) \\
		t & \mbox{otherwise}
	\end{array}
\right.
$$

Suppose we have a $\Sigma$-algebra $(X, c)$, a term $t \in T$, and a position $p \in Pos(t)$. We would like to find a function $\alpha_p^t: X \rightarrow X$ which satisfies
$$
\begin{tikzcd}[row sep=huge,column sep=huge]
T
\arrow[r, "Sub_p^t"]
\arrow[d, "\cata(c)"]
& T
\arrow[d, "\cata(c)"]
\\
X
\arrow[r, "\alpha_p^t"]
& X
\end{tikzcd}
$$
and do so by induction on $l(p)$. If $l(p) = 0$ then $p = \epsilon$ and $Sub_p^t = id_T$, so $\alpha_p^t = id_X$. If $l(p) \geq 1$ then $p = kq$ for some $k \in \N$ and position $q$. Since $p \in Pos(t)$ we must have $d(t) \geq 1$, so $t = f(t_1,...,t_n)$ for some $f \in \Sigma_n$, where $n \geq k$. Now $l(q) = l(p)-1$, so by induction there is some $\alpha_q^{t_k}: X \rightarrow X$ which satisfies
$$
\begin{tikzcd}[row sep=huge,column sep=huge]
T
\arrow[r, "Sub_q^{t_k}"]
\arrow[d, "\cata(c)"]
& T
\arrow[d, "\cata(c)"]
\\
X
\arrow[r, "\alpha_q^{t_k}"]
& X
\end{tikzcd}
$$
So we define $\alpha_p^t: X \rightarrow X$ by
$$
x \mapsto c(f)(\cata(c)(t_1),...,\cata(c)(t_{k-1}),\alpha_q^{t_k}(x),\cata(c)(t_{k+1}),...,\cata(c)(t_n)).
$$
We can then verify that, if $s \in T$, we have
\begin{equation*}
    \begin{split}
        \alpha_p^t\circ \cata(c)(s)
        &= c(f)(\cata(c)(t_1),...,\alpha_q^{t_k}(\cata(c)(s)),...,\cata(c)(t_n)) \\
        &= c(f)(\cata(c)(t_1),...,\cata(c)(t_k[s]_q),...,\cata(c)(t_n)) \\
        &= \cata(c)(f(t_1,...,t_k[s]_q,...,t_n)) \\
        &= \cata(c)(t[s]_{kq}) \\
        &= \cata(c)\circ Sub_p^t(s),
    \end{split}
\end{equation*}
so $\alpha_p^t$ satisfies $\alpha_p^t\circ \cata(c) = \cata(c)\circ Sub_p^t$.

We can now expand the diagram given above with these additional pieces, where $t_0$ is the initial term (which lives in $T_p^l$).
$$
\begin{tikzcd}[row sep=huge,column sep=huge]
T_p^l
\arrow[rrr, "R_p^n"]
\arrow[d, "r_p"]
&&& T_p^l
\arrow[d, "r_p"]
&
\\
T_\epsilon^l
\arrow[r, "\Psi_{v(l)}"]
& T^{v(l)}
\arrow[r, "H^n"]
\arrow[d, "F^{v(l)}(\cata(c))"]
& T^{v(l)}
\arrow[d, "F^{v(l)}(\cata(c))"]
\arrow[r, "\Phi_{v(l)}"]
& T_\epsilon^l
\arrow[d, "\cata(c)"]
\arrow[r, "Sub_p^{t_0}"]
& T_p^l
\arrow[d, "\cata(c)"]
\\
& X^{v(l)}
\arrow[r, "G^n"]
& X^{v(l)}
\arrow[r, "C_{v(l)}"]
& X
\arrow[r, "\alpha_p^{t_0}"]
& X
\end{tikzcd}
$$
We notice that, since $r = \tau(l)$, if $t|_p = \sigma(l)$ then $R_p(t) = t[\sigma\tau(l)]_p$ and so $R_p^n(t) = t[\sigma\tau^n(l)]_p$. So given the initial term $t_0\in T_p^l$ where $t_0|_p = \sigma(l)$, we have
\begin{equation*}
    \begin{split}
        Sub_p^{t_0} \circ r_p \circ R_p^n(t_0)
        &= Sub_p^{t_0}\circ r_p (t_0[\sigma\tau^n(l)]_p) \\
        &= Sub_p^{t_0}(\sigma\tau^n(l)) \\
        &= t_0[\sigma\tau^n(l)]_p \\
        &= R_p^n(t_0).
    \end{split}
\end{equation*}
If we now let
$$
x_0 := F^{v(l)}(\cata(c)) \circ \Psi_{v(l)} (t_0|_p),
$$
be our initial element of $X^{v(l)}$, we can conclude that
\begin{equation*}
    \begin{split}
        \cata(c) \circ R_p^n(t_0)
        &= \cata(c) \circ Sub_p^{t_0} \circ r_p \circ R_p^n(t_0) \\
        &= \alpha_p^{t_0} \circ C_{v(l)} \circ G^n \circ F^{v(l)} (\cata(c)) \circ \Psi_{v(l)} \circ r_p (t_0) \\
        &= (\alpha_p^{t_0} \circ C_{v(l)}) \circ G^n (x_0).
    \end{split}
\end{equation*}
So the rewriting model is again a cartesian dynamical system, where the hidden dynamics are still given by $G$, but where $\alpha_p^{t_0} \circ C_{v(l)}$ is the new \emph{output} map $X^{v(l)}\rightarrow X$.

We now prove the converse: that the cartesian dynamical systems embed in the rewriting models, in the sense that composition with the construction above will fix the dynamical system. Suppose we have a cartesian dynamical system on $X$ of depth $d$, so there are maps $G: X^d \rightarrow X^d$ and $f: X^d \rightarrow X$. We would like to find an equivalent rewriting model, i.e. a section (left-inverse) of the projection above. Define a signature $\Sigma$ via
$$
\Sigma_0 = \{a_1,...,a_d\} \qquad \Sigma_1 = \{\iota\} \qquad \Sigma_d = \{\sigma_0,...,\sigma_d\}
$$
and $\Sigma_n = \varnothing$ for all other $n$. Now $F(X) = d + X + (d+1)X^d$. We pick a variable set $V = \{v_1,...,v_d\}$ and a $\Sigma$-identity
$$
l = \sigma_0(v_1,...,v_d)
\qquad
r = \sigma_0(\sigma_1(v_1,...,v_d), ..., \sigma_d(v_1,...,v_d))
$$
which induces a natural transformation $\eta: F \rightarrow F^2$ where
$$
\sigma_0(t_1,...,t_d) \mapsto \sigma_0(\sigma_1(t_1,...,t_d), ..., \sigma_d(t_1,...,t_d))
$$
$$
a_i \mapsto \iota^2(a_i)
$$
$$
\sigma_i(t_1,...,t_d) \mapsto \iota(\sigma_i(a_1,...,a_d))
$$
for $i = 1,...,d$. We also induce a catamorphism by the assignments
$$
\sigma_0 \mapsto f \qquad \sigma_i \mapsto \pi_i \circ G \qquad \iota \mapsto id_X \qquad a_i \mapsto x_i
$$
for $i = 1,...,d$, where $x_0 = (x_1,...,x_d) \in X^d$ is the initial internal state of the dynamical system, and $\pi_i:X^d \rightarrow X$ is projection onto the $i^{th}$ component. We can represent this in the following diagram.
$$
\begin{tikzcd}[row sep=huge,column sep=large]
T
\arrow[r, "\psi"]
\arrow[rrr, "R_\epsilon^n", bend left=20]
& d + T + (d+1)T^d
\arrow[r, "H^n"]
\arrow[d, "F(\cata(c))"]
& d + T + (d+1)T^d
\arrow[d, "F(\cata(c))"]
\arrow[r, "\phi"]
& T
\arrow[d, "\cata(c)"]
\\
& d + X + (d+1)X^d
\arrow[r, "G^n"]
& d + X + (d+1)X^d
\arrow[r, "c"]
& X
\end{tikzcd}
$$
If we choose the initial term $t_0\in T$ as $t_0 = \sigma_0(a_1,...,a_d)$ then
\begin{equation*}
    \begin{split}
        F(\cata(c)) \circ \psi (t_0) &= F(\cata(c)) (a_1,...,a_d) \\ 
        &= (\cata(c)(a_1),...,\cata(c)(a_d)) \\
        &= (x_1,...,x_d) \\
        &= x_0.
    \end{split}
\end{equation*}
and so $\cata(c)\circ R_\epsilon^n(t_0) = f \circ G^n (x)$ for all $n$. So any cartesian dynamical system embeds as a rewriting model at position $\epsilon$.
\end{proof}

We can interpret this result as saying that we can use the purely algebraic language of terms and rewriting to describe the structure of a dynamical system, but that the structure contained in the rewrite rule is strictly richer than the dynamical system projection. This allows us to identify how the algebraic structure evolves and whether it has any interesting properties. The converse is also true: given a rewriting model for a sequence in an arbitrary type $X$, we can now understand the model to be the appropriate generalisation of dynamical systems to that type.

\section{Compositionality of Rewriting Models}

Given that we now have an algebraic language for the temporal structure of these models, we would also like to know in what sense it is compositional. In other words, how do we know what properties of a dynamical system are propagated over time? In order to address this we identify compositionality with \emph{functoriality}, meaning that the properties we want to be preserved under composition define a category $C$ which embeds via a forgetful functor $U$ into $\Set$.

Examples of this can be constraints on the sorts of sets involved, which could be topological spaces, manifolds or Euclidean spaces. The functions between them can then be further reduced to be continuous, smooth, $k$-smooth, or linear. Another important class of categories are the $G$-equivariant sets, which are the main object of \emph{geometric deep learning}, where models are prescribed to be equivariant under the actions or representations of various groups $G$. Convolutional neural networks, for example, are maps in the category of $\R^2$-equivariant Euclidean spaces with piecewise linear maps. This is achieved via compositionality, since each layer can be shown to be $\R^2$-equivariant, and so their composition must be too.

We will prove the following meta-theorem, which states that, if the sets and maps that define the dynamical system are in the image of the functor $U$, then, under certain assumptions on $C$, the dynamical system as a whole will lift to one in $C$.

\begin{theorem}[Compositionality] \label{compositionality}
Let $C$ be a category with products and coproducts, and $U:C \rightarrow \Set$ be a functor with natural isomorphisms
$$
U(X+Y) \cong U(X)+U(Y)
\qquad
U(X\times Y) \cong U(X)\times U(Y).
$$
Let $X\in C$ and $c_\sigma:X^{|\sigma|}\rightarrow X$ in $C$ for each $\sigma$ in some signature $\Sigma$, so that we can define $c:F(U(X)) \rightarrow U(X)$ with each $U(c_\sigma)$. Then a rewriting model given by some $\Sigma$-identity $(l,r)$ and $c$ is equivalent to a cartesian dynamical system whose sets and maps are all in $C$.
\end{theorem}
\begin{proof}
The statement is well-defined because $U$ preserves products, so if $c_\sigma:X^{|\sigma|}\rightarrow X$ in $C$ then $U(c_\sigma):U(X)^{|\sigma|}\rightarrow U(X)$ in $\Set$. The fact that $U$ also preserves coproducts means we can define a lift $\hat{F}:C \rightarrow C$ such that $U\hat{F} = FU$, where $\hat{F}$ copies the same polynomial form as $F$. This is possible since the polynomials induced by a rewriting model have finite exponents, and we can define $X^n$ inductively in $C$.

Suppose the $\Sigma$-identity $(l,r)$ induces a natural transformation $\eta:F^{L(l)} \rightarrow F^{L(r)}$. We would like to construct a lift of $\eta_X$ to a map $\hat{\eta}:\hat{F}^{L(l)}(X)\rightarrow\hat{F}^{L(r)}(X)$ in $C$ (note that we just require a map: this $\hat{\eta}$ will not be a natural transformation). If we denote by $Y^k$ the functor $Y \mapsto Y^k$ then we can decompose $\eta$ into natural transformations $Y^n \rightarrow Y^m$. Using the Yoneda lemma we identify each of these with maps $u:m \rightarrow n$ and can express the action of $\eta_X$ as $x\mapsto x \circ u$. We can then define $\hat{\eta}$ on each cofactor $X^n$ by sending $x\mapsto x \circ u \in X^m$. It is clear that $U(\hat{\eta}) = \eta_{U(X)}$ by construction.

Finally, we note that $c:F(U(X)) \rightarrow U(X)$ is constructed such that $c = U(\hat{c})$ where $\hat{c}:\hat{F}(X) \rightarrow X$ is given by the maps $c_\sigma:X^{|\sigma|}\rightarrow X$ on each cofactor of $\hat{F}(X)$. It then also follows that $F(c) = U(\hat{F}(\hat{c}))$. So all the maps involved in the construction of the rewriting model are in the image of $U$, and we can build a lifted model in $C$ whose image under $U$ returns the original model. The fact that $U$ preserves coproducts means that the lifted model in $C$ is also a cartesian dynamical system.
\end{proof}

In the context of functional programming or machine learning, the categories of interest $C$ are all modelled in $\Set$, meaning that a forgetful functor $U : C \rightarrow \Set$ is always defined, and will usually have a left adjoint \emph{free} functor $\Set \rightarrow C$. Since right adjoints preserve limits, the product in $C$ will agree with the one in set and so the product condition in Theorem \ref{compositionality} is immediately satisfied. The subtlety is in the coproduct condition, which is generally satisfied when $C$ is topological, like smooth manifolds or Euclidean spaces. In algebraic categories like vector spaces or groups, the condition will break down, as the coproduct has to identify the identity elements of its cofactors, and so we cannot define maps separately on each of them.

It follows that, even though we have defined rewriting models in the minimal setting of $\Set$, they are actually a universal construction in all categories with the appropriate product and coproduct structures. This includes the topological categories with continuous or smooth maps, and the important group-equivariant topological categories of interest in geometric deep learning. So in the context of machine learning, this means that rewriting models are compositional. Specifically, if the sets and functions involved all have a desirable property (i.e. form a category), then the model as a whole also has this property.

Finally, we observe that we can make the totally recursive structure of a rewriting model explicit as a composition of catamorphisms over $\N$. We can describe the iterated application of $R_p$ to the term $t_0$ as the catamorphism
$$
\begin{tikzcd}[row sep=huge,column sep=huge]
1+\N
\arrow[d, "0+succ"]
\arrow[r, "1 + \cata(t_0 + R)"]
& 1+T_p^l
\arrow[d, "t_0 + R_p"]
\\
\N
\arrow[r, "\cata(t_0 + R_p)"]
& T_p^l
\end{tikzcd}
$$
so that the model is the composition
$$
\cata(c) \circ \cata(t_0 + R_p) : \N \rightarrow T_p^l \rightarrow X.
$$

\backmatter

\bmhead{Acknowledgements}

We would like to acknowledge the generous financial support and opportunities for collaboration that have been given to this work by Hylomorph Solutions. Jeffrey Giansiracusa was supported by EPSRC grant EP/R018472/1 through the Centre for TDA.

\section*{Declarations}

\subsection*{Ethics approval} Not applicable.

\subsection*{Consent to participate} Not applicable.

\subsection*{Consent for publication} All authors consent to the publication of this work.

\subsection*{Availability of data and materials} Not applicable.

\subsection*{Competing interests} 
Financial interests: Iolo Jones is a consultant to Hylomorph Solutions. Jerry Swan receives a salary from Hylomorph Solutions where he is Director of Research. Jeffrey Giansiracusa is funded by EPSRC grant EP/R018472/1 through the Centre for TDA. Non-financial interests: None.

\subsection*{Funding} Jeffrey Giansiracusa was supported by EPSRC grant EP/R018472/1 through the Centre for TDA.

\subsection*{Authors' contributions}
Iolo Jones and Jerry Swan conceived of the rewriting model framework. The main manuscript was written by Iolo Jones with the supervision of Jeffrey Giansiracusa, and the introduction was written by Iolo Jones and Jerry Swan. All authors reviewed the manuscript.

\printbibliography


\end{document}